%% file: main.tex
\documentclass[10pt]{amsart}
\title{Poincare Series and Miraculous Duality}
\author{Kevin Lin}

\usepackage[margin=1.5in]{geometry}
\usepackage{etoolbox}
\usepackage[nameinlink]{cleveref}

\usepackage{baskervald}
\usepackage[T1]{fontenc}

\usepackage[mathscr]{euscript}

\usepackage{amssymb}
\usepackage{tikz-cd}

\usepackage[style=alphabetic,doi=false,isbn=false,url=false,eprint=false,sorting=nyt]{biblatex}
\numberwithin{equation}{subsubsection}

\newcommand\newkevintheorem[3][]{%
    \newtheorem{#2}[kevinthm]{#3}%
    \csletcs{old#2}{#2}%
    \expandafter\renewcommand\csname #2\endcsname{\stepcounter{subsubsection}\csname old#2\endcsname}%
}

\newtheorem{innercustomgeneric}{\customgenericname}
\providecommand{\customgenericname}{}
\newcommand{\newcustomtheorem}[2]{%
  \newenvironment{#1}[1]
  {%
   \renewcommand\customgenericname{#2}%
   \renewcommand\theinnercustomgeneric{\kern-0.3em ##1}%
   \innercustomgeneric
  }
  {\endinnercustomgeneric}
}

\newcustomtheorem{specialtheorem}{}
\crefname{specialtheorem}{}{}

\crefname{appsec}{Appendix}{Appendices}

\newtheoremstyle{kevintheoremstyle}
{}
{}
{}
{0pt}
{\bfseries}
{.}
{4pt}
{\thmnumber{#2}\thmname{ #1}\thmnote{ (#3)}}

\theoremstyle{kevintheoremstyle}

\newkevintheorem{theorem}{Theorem}
\crefname{theorem}{Theorem}{Theorems}
\newkevintheorem{lemma}{Lemma}
\crefname{lemma}{Lemma}{Lemmas}
\newkevintheorem{proposition}{Proposition}
\crefname{proposition}{Proposition}{Propositions}
\newkevintheorem{corollary}{Corollary}
\crefname{corollary}{Corollary}{Corollaries}
\newkevintheorem{construction}{Construction}
\crefname{construction}{Construction}{Constructions}

\newkevintheorem{definition}{Definition}
\crefname{definition}{Definition}{Definitions}

\input{commands}

\begin{document}

\maketitle

\begin{abstract}
In the setting of global geometric Langlands, we show that miraculous duality on the stack of principal bundles on a curve intertwines the functor of Poincare series with the dual functor to Whittaker coefficients. We construct, for arbitrary parabolic subgroups, Jacquet functors controlling constant terms of Poincare series and Whittaker coefficients of Eisenstein series in local terms.
\end{abstract}

\input{introduction}
\input{construction}
\input{cofiber}
\input{Jacquet}
\input{ULA}
\input{coJacquet}
\input{Ran}
\input{Ran_upgrades}
\input{3_morphism}

\printbibliography

\end{document}

%% file: commands.tex
\newcommand{\im}{\mathrm{im}}

\newcommand{\Maps}{\mathrm{Maps}}

\newcommand{\Spec}{\mathrm{Spec}}

\newcommand{\id}{\mathbf{1}}
\newcommand{\fib}{\mathrm{fib}}
\newcommand{\cofib}{\mathrm{cofib}}

\newcommand{\mon}{\hookrightarrow}
\newcommand{\epi}{\twoheadrightarrow}

\newcommand{\ev}{\mathrm{ev}}
\newcommand{\AJ}{\mathrm{AJ}}
\newcommand{\AS}{\mathrm{AS}}
\newcommand{\Av}{\mathrm{Av}}
\newcommand{\ext}{\mathrm{ext}}
\newcommand{\oblv}{\mathrm{oblv}}
\newcommand{\pr}{\mathrm{pr}}

\renewcommand{\AA}{\mathbf{A}}

\newcommand{\BB}{\mathbf{B}}

\newcommand{\DD}{\mathbf{D}}
\newcommand{\Deltac}{\overline{\Delta}}
\newcommand{\EE}{\mathscr{E}}
\newcommand{\FF}{\mathscr{F}}
\newcommand{\FFo}{\overset{\circ}{\FF}{}}
\newcommand{\GG}{\mathscr{G}}
\newcommand{\LG}{{\check{G}}}
\newcommand{\Lg}{{\check{\fg}}}
\newcommand{\fg}{\mathfrak{g}}
\newcommand{\Ga}{\mathbf{G}_{\mathrm{a}}}
\newcommand{\GL}{\mathrm{GL}}
\newcommand{\Gm}{\mathbf{G}_m}
\newcommand{\rH}{\mathrm{H}}
\newcommand{\HH}{\mathscr{H}}
\newcommand{\HHt}{\widetilde{\HH}}

\newcommand{\LL}{\mathscr{L}}
\newcommand{\GLL}{\mathbf{L}}

\newcommand{\OO}{\mathscr{O}}

\newcommand{\bfP}{\mathbf{P}}
\newcommand{\fp}{\mathfrak{p}}

\renewcommand{\SS}{\mathscr{S}}
\newcommand{\SL}{\mathrm{SL}}

\newcommand{\UU}{\mathscr{U}}
\newcommand{\fu}{\mathfrak{u}}
\newcommand{\WW}{\mathscr{W}}
\newcommand{\YY}{\mathscr{Y}}
\newcommand{\ZZ}{\mathscr{Z}}
\newcommand{\ZZt}{{\widetilde{\ZZ}}}

\newcommand{\Bc}{\overline{B}}
\newcommand{\Nc}{\overline{N}}
\newcommand{\Pt}{\widetilde{P}}

\newcommand{\qt}{\widetilde{q}}

\newcommand{\aph}{{'}}

\newcommand{\rhov}{{\check{\rho}}}

\newcommand{\ab}{\mathrm{ab}}
\newcommand{\ad}{\mathrm{ad}}
\newcommand{\cis}{\mathrm{cis}}
\newcommand{\co}{\mathrm{co}}
\newcommand{\cusp}{\mathrm{cusp}}
\newcommand{\disj}{\mathrm{disj}}

\newcommand{\gen}{\mathrm{gen}}
\newcommand{\glob}{\mathrm{glob}}
\newcommand{\good}{\mathrm{good}}
\newcommand{\gr}{\mathrm{gr}}

\renewcommand{\neg}{\mathrm{neg}}
\newcommand{\Nilp}{\mathscr{N}}

\newcommand{\pol}{\mathrm{pol}}
\newcommand{\vpol}{{\mathrm{v.}\pol}}
\newcommand{\pos}{\mathrm{pos}}

\newcommand{\ren}{\mathrm{ren}}
\newcommand{\rev}{\mathrm{rev}}
\newcommand{\spec}{\mathrm{spec}}
\newcommand{\trans}{\mathrm{tr}}

\newcommand{\vac}{\mathrm{vac}}

\newcommand{\Bun}{\mathrm{Bun}}
\newcommand{\Bunc}{\overline{\mathrm{Bun}}}
\newcommand{\Bunt}{\widetilde{\mathrm{Bun}}}
\newcommand{\VBun}{\mathrm{VinBun}}
\newcommand{\Conf}{\mathrm{Conf}}
\newcommand{\Gr}{\mathrm{Gr}}
\newcommand{\Grc}{\overline{\Gr}}
\newcommand{\Grt}{\widetilde{\Gr}}

\newcommand{\Mod}{\mathrm{Mod}}
\newcommand{\LS}{\mathrm{LS}}
\newcommand{\Ran}{\mathrm{Ran}}
\newcommand{\Sym}{\mathrm{Sym}}
\newcommand{\Vin}{\mathrm{Vin}}

\newcommand{\DMod}{\mathrm{D}}
\newcommand{\ICoh}{\mathrm{IndCoh}}

\newcommand{\QCoh}{\mathrm{QCoh}}
\newcommand{\Sh}{\mathrm{Sh}}
\newcommand{\Sing}{\mathrm{SS}}

\newcommand{\CS}{\mathrm{CS}}
\newcommand{\Rep}{\mathrm{Rep}}
\newcommand{\Sat}{\mathrm{Sat}}
\newcommand{\Sph}{\mathrm{Sph}}
\newcommand{\Vect}{\mathrm{Vect}}
\newcommand{\Whit}{\mathrm{Whit}}

\newcommand{\coeff}{\mathrm{coeff}}
\newcommand{\CT}{\mathrm{CT}}
\newcommand{\DL}{\mathrm{DL}}
\newcommand{\Eis}{\mathrm{Eis}}
\newcommand{\Jac}{\mathrm{Jac}}
\newcommand{\Jact}{\widetilde{\Jac}}
\newcommand{\PsId}{\mathrm{PsId}}
\newcommand{\Poinc}{\mathrm{Poinc}}

\newcommand{\sprd}{\mathrm{sprd}}

\newcommand{\Disc}{\mathscr{D}}
\newcommand{\Disco}{\overset{\circ}{\Disc}}
\newcommand{\pt}{\mathrm{pt}}

\newcommand{\bs}{\backslash}
\newcommand{\bss}{\bs\!\!\bs}
\newcommand{\git}{/ \!\! /}

\newcommand{\simto}{\overset{\sim}{\rightarrow}}

\newcommand{\timest}{\widetilde{\times}}

\newcommand{\shimes}{\overset{!}{\otimes}}
\newcommand{\stimes}{\overset{*}{\otimes}}

\newcommand{\shoximes}{\overset{!}{\boxtimes}}
\newcommand{\stoximes}{\overset{*}{\boxtimes}}

\newcommand{\LHS}{\mathrm{LHS}}
\newcommand{\RHS}{\mathrm{RHS}}

\newcommand{\tightoverset}[2]{%
  \mathop{#2}\limits^{\vbox to -.5ex{\kern-0.75ex\hbox{$#1$}\vss}}}

\newcommand{\hl}{{\tightoverset{\leftharpoonup}{h}}}
\newcommand{\hr}{{\tightoverset{\rightharpoonup}{h}}}

%% file: introduction.tex
\section{Introduction}

\subsection{Statement of the theorem}

\subsubsection{} We fix a smooth proper curve $X$ over an algebraically closed field of characteristic zero, as well as a square root $\omega_X^{1/2}$ of the canonical line bundle.

\subsubsection{} The geometric Langlands conjecture \cite{ss_indcoh} predicts an equivalence of categories
    \[\GLL_G : \DMod(\Bun_G) \simto \ICoh_\Nilp(\LS_\LG).\]
Let us review some of the compatibilities expected from this functor.

\subsubsection{Duality} The category appearing on the spectral side of the Langlands equivalence has a self-duality given by Serre duality
    \[
    \DD : \ICoh_\Nilp(\LS_\LG) \simto \ICoh_\Nilp(\LS_\LG)^\vee
    \qquad
    \FF \mapsto \Gamma (X, \FF \shimes -)
    \]
The work \cite{strange} constructs a self duality
    \[
    \PsId_! : \DMod(\Bun_G)^\vee \simto \DMod(\Bun_G)
    \]
which is expected to be the automorphic counterpart of Serre duality in the sense that the following diagram commutes up to a Cartan involution and cohomological shift.
\begin{equation}
\label{eqn:Serre_PsId_no_shift}
\begin{tikzcd}[
  column sep={9em,between origins},
  row sep={4em,between origins}
]
  \DMod(\Bun_G) \arrow[r, "\GLL_G"] \arrow[dr, phantom, "\tau \cdot {[?]}"] & \ICoh_\Nilp(\LS_\LG) \arrow[d, "\DD"] \\
  \DMod(\Bun_G)^\vee \arrow[u, "\PsId_!"] & \ICoh_\Nilp(\LS_\LG)^\vee \arrow[l, "\GLL_G^\vee"]
\end{tikzcd}
\end{equation}
By definition, $\PsId_!$ is given by the object
    \[
    \Delta_! \cdot k_{\Bun_G} \in
    \DMod(\Bun_G \times \Bun_G)
    \overset{\sim}{\leftarrow} \DMod(\Bun_G) \otimes \DMod(\Bun_G),
    \]
where $\Delta : \Bun_G \to \Bun_G \times \Bun_G$ is the diagonal.

\subsubsection{Whittaker coefficients} Let $\rhov$ be the half-sum of the positive coroots. By using the chosen square root of the canonical bundle, we obtain a point of $\Bun_T$ denoted $\omega^\rhov$ and define
    \[
    \Bun_N = \{\omega^\rhov\} \underset{\Bun_T}{\times} \Bun_B.
    \]
This stack is equipped with canonical maps
    \[
    \begin{tikzcd}[
  column sep={3em,between origins},
  row sep={3em,between origins}
]
      & \Bun_N \arrow[dl, "\psi" above left] \arrow[dr, "\pi"] & \\
      \Ga & & \Bun_G
    \end{tikzcd}
    \]
Now let $\AS$ denote the exponential D-module on $\Ga$, cohomologically shifted to be a character sheaf with respect to $!$-convolution. We write
    \[\Vect \simto \Whit(\Bun_N) \overset{\oblv}{\longrightarrow} \DMod(\Bun_N)\]
for the full subcategory generated by the object $\WW_\vac = \ev^! \cdot \AS$. The inclusion admits a right adjoint
    \[
    \DMod(\Bun_N) \overset{\Av_*}{\longrightarrow} \Whit(\Bun_N) \simeq \Vect,
    \]
and the functor of Whittaker coefficients is
    \[
    \coeff_* = \Av_* \cdot \pi^!.
    \]
It is supposed to match the functor of pushforward along $f : \LS_\LG \to \pt$.
\begin{equation}
\label{eqn:L_G_coeff_Poinc}
\begin{tikzcd}[
  column sep={8em,between origins},
  row sep={4em,between origins}
]
  \DMod(\Bun_G) \arrow[d, "\coeff_*" left] \arrow[r, "\GLL_G"] & \ICoh_\Nilp(\LS_\LG) \arrow[d, "f_*"] \\
  \Whit(\Bun_N) \arrow[r, "\sim"] & \QCoh(\pt)
\end{tikzcd}
\qquad
\begin{tikzcd}[
  column sep={8em,between origins},
  row sep={4em,between origins}
]
    \Whit(\Bun_N) \arrow[r, "\sim"] \arrow[d, "\Poinc_!" left] & \QCoh(\pt) \arrow[d, "f^*"] \\ 
    \DMod(\Bun_G) \arrow[r, "\GLL_G"] & \ICoh_\Nilp(\LS_\LG)
\end{tikzcd}
\end{equation}
Passing to left adjoints gives the diagram on the right, where
    \[\Poinc_! = \pi_! \cdot \oblv.\]

\subsubsection{} It is a general fact that Serre duality exchanges $f_*$ with $f^!$. However, a special feature of $\LS_\LG$ is that there is a canonical isomorphism
    \[\omega_{\LS_\LG} \simeq \OO_{\LS_\LG} \, [2 \delta_G]\]
where
    \[2 \delta_G = (2g - 2) \cdot \dim(\LG) = 2 \dim \Bun_G.\]
Therefore the dual functor to $f_*$ identifies with a shift of $f^*$.
\[
\begin{tikzcd}[
  column sep={10em,between origins},
  row sep={4em,between origins}
]
    \QCoh(\pt) \arrow[r, "\sim"] \arrow[d, "f^*{[2 \delta_G]}" left] & \QCoh(\pt)^\vee \arrow[d, "(f_*)^\vee"] \\ 
    \ICoh_\Nilp(\LS_\LG) \arrow[r, "\DD"] & \ICoh_\Nilp(\LS_\LG)^\vee
\end{tikzcd}
\]
Let us translate this statement to the automophic side. We will use the Verdier duality pairing
    \[
    \Whit_{-\psi}(\Bun_N) \otimes \Whit_\psi(\Bun_N) \to \Vect
    \qquad
    \langle \WW_1, \WW_2 \rangle = \Gamma(\Bun_N, \WW_1 \shimes \WW_2).
    \]
The resulting statement is our main theorem:

\begin{theorem}
There is a canonical commutative diagram
\[
\label{thm:vac}
\begin{tikzcd}[
  column sep={10em,between origins},
  row sep={4em,between origins}
]
    \Whit_{-\psi}(\Bun_N)^\vee \arrow[d, "\coeff_*^\vee {[2 \delta_G]}" left] \arrow[r] & \Whit_\psi(\Bun_N) \arrow[d, "\Poinc_!"] \\
    \DMod(\Bun_G)^\vee \arrow[r, "\PsId_!"] & \DMod(\Bun_G)
\end{tikzcd}
\]
\end{theorem}

\subsubsection{} As we have just explained, Theorem \ref{thm:vac} is predicted by the geometric Langlands conjecture. However, we expect that it will play an essential r\^ole in the construction of the Langlands equivalence $\GLL_G$ and its attendant compatibilities. 

\subsubsection{Cohomological shifts} Due the ambiguity in (\ref{eqn:Serre_PsId_no_shift}), the above discussion does not actually predict the cohomological shift appearing in our theorem. (One can guess this shift by considering the compatibility expressed by Proposition \ref{proposition:poinc_identification}, which only involves the automorphic side of the Langlands correspondence.) Let us use the theorem statement to determine the shift that should appear in (\ref{eqn:Serre_PsId_no_shift}). The Verdier duality on $\Whit(\Bun_N)$ is related to Serre duality as follows:
\[
\begin{tikzcd}[
  column sep={9em,between origins},
  row sep={4em,between origins}
]
  \Whit_\psi(\Bun_N) \arrow[r, "\sim"] \arrow[d, "\DD{[-2 \delta_N]}" left] & \QCoh(\pt) \arrow[d, "\DD"] \\
  \Whit_{-\psi}(\Bun_N)^\vee \arrow[r, "\sim"] & \QCoh(\pt)^\vee
\end{tikzcd}
\]
The shift appearing here is actually to be expected: it accounts for the difference between Verdier duality and the more natural local Whittaker duality. In any case we obtain:
\[
\begin{tikzcd}[
  column sep={9em,between origins},
  row sep={4em,between origins}
]
  \DMod(\Bun_G) \arrow[r, "\GLL_G"] \arrow[dr, phantom, "\tau"] & \ICoh_\Nilp(\LS_\LG) \arrow[d, "\DD {[2\delta_G]}"] \\
  \DMod(\Bun_G)^\vee \arrow[u, "\PsId_! {[2 \delta_G + 2 \delta_N]}"] & \ICoh_\Nilp(\LS_\LG)^\vee \arrow[l, "\GLL_G^\vee"]
\end{tikzcd}
\]

\subsection{Upgraded version}

\subsubsection{} Our main theorem has a $\Sph_{G,\Ran}$-linear upgrade involving the entire factorization category $\Whit(G)_\Ran$. Let us recall some basic facts about Whittaker categories and Hecke actions before giving the statement.

\subsubsection{Global Whittaker model} Let $x$ be a finite set of points of $X$. We will use the global model for the Whittaker category
    \[
    \Whit(G)_x \simeq \Whit\left(\Bunc_{N,x}^\pol\right) \mon \DMod\left(\Bunc_{N,x}^\pol \right).
    \]
It is naturally a module over the spherical Hecke category $\Sph_{G,x}$, equipped with a pair of $\Sph_{G,x}$-linear adjoint functors
    \[
    \Poinc_! : \Whit(G)_x \to \DMod(\Bun_G)
    \qquad
    \coeff_* : \DMod(\Bun_G) \to \Whit(G)_x.
    \]
Here, as before, $\Poinc_!$ is simply $!$-pushforward of the underlying sheaf to $\Bun_G$. The Verdier pairing
    \[
    \Whit_{-\psi}(G)_x \otimes \Whit_\psi(G)_x \to \Vect
    \]
is a perfect pairing.

\subsubsection{Left and right modules} To formulate the sense in which our dualities are $\Sph_G$-linear, recall that every convolution stack comes with a canonical involution
    \[
    \sigma : X \underset{S}{\times} X \to X \underset{S}{\times} X
    \]
given by swapping the two factors. For example, this involution on $G \simeq \pt \underset{\BB G}{\times} \pt$ is group inversion. We obtain an involutive anti-automorphism
    \[\sigma : \Sph_{G,x} \simto \Sph_{G,x}^\rev.\]
Now let us view $\DMod(\Bun_G)$ as a right $\Sph_{G,x}$-module. Then $\DMod(\Bun_G)^\vee$ is naturally a left $\Sph_{G,x}$-module, but the involution $\sigma$ allows us to identify left and right modules for $\Sph_{G,x}$. Therefore, it makes sense to say that the functor
    \[\PsId_! : \DMod(\Bun_G)^\vee \to \DMod(\Bun_G)\]
is canonically $\Sph_{G,x}$-linear. This is easy to check using the ULA property of Hecke kernels and the fact that $\sigma$ interchanges the geometrically defined left and right $\Sph_{G,x}$ actions on $\DMod(\Bun_G)$. It is tautological that the duality
    \[\Whit_{-\psi}(G)_x^\vee \simeq \Whit_\psi(G)_x\]
is $\Sph_{G,x}$-linear.

\subsubsection{} Allowing the set of points $x$ to vary in families and writing 
    \[\DMod(\Bun_G)_\Ran = \DMod(\Bun_G \times \Ran),\]
we obtain the following variant of the main theorem:

\begin{theorem}
\label{thm:Ran}
There is a canonical $\Sph_{G,\Ran}$-linear commutative diagram
\[
\begin{tikzcd}[
  column sep={10em,between origins},
  row sep={4em,between origins}
]
    \Whit_{-\psi}(G)_\Ran^\vee \arrow[d, "\coeff_*^\vee {[2 \delta_G]}" left] \arrow[r] & \Whit_\psi(G)_\Ran \arrow[d, "\Poinc_!"] \\
    \DMod(\Bun_G)_\Ran^\vee \arrow[r, "\PsId_!"] & \DMod(\Bun_G)_\Ran
\end{tikzcd}
\]
\end{theorem}

\subsubsection{} Let us explain why this theorem is to be expected from the point of view of Langlands duality. In our derivation of the main theorem for the Whittaker vacuum, we actually ignored the Cartan involution without comment. This was acceptable because the Cartan involution sends $\OO_{\LS_\LG}$ to itself. But in what follows we will have to track this involution more closely.

\subsubsection{} Recall that geometric Satake defines a monoidal functor
    \[\Sat_x : \Rep(\LG)_x \to \Sph_{G,x}.\]
The Langlands functor is supposed to intertwine the Hecke action of $\Sph_{G,x}$ with the action of $\Rep(\LG)_x$ on $\ICoh_\Nilp(\LS_\LG)$ induced by the map
    \[\ev_x : \LS_\LG(X) \to \LS_\LG(\mathscr{D}_x) \simeq (\BB \LG)_x\]
and the action of $\QCoh$ on $\ICoh$. Therefore, viewing $\WW_\vac$ as an object of $\Whit(G)_x$ and acting on (\ref{eqn:L_G_coeff_Poinc}), we obtain:
\[
\begin{tikzcd}[
  column sep={8em,between origins},
  row sep={4em,between origins}
]
    \Whit(G)_x \arrow[d, "\Poinc_!"] & \Rep(\LG)_x \arrow[d, "\ev_x^*"] \arrow[l, "\CS" above] \\ 
    \DMod(\Bun_G) \arrow[r, "\GLL_G" below] & \ICoh_\Nilp(\LS_\LG)
\end{tikzcd}
\qquad
\begin{tikzcd}[
  column sep={8em,between origins},
  row sep={4em,between origins}
]
  \DMod(\Bun_G) \arrow[d, "\coeff_*"] \arrow[r, "\GLL_G"] & \ICoh_\Nilp(\LS_\LG) \arrow[d, "\ev_{x,*}"] \\
  \Whit(G)_x & \Rep(\LG)_x \arrow[l, "\CS"]
\end{tikzcd}
\]

\subsubsection{Cartan involution} Let us recall the behavior of the involution $\sigma$ with respect to geometric Satake. To reduce confusion, we will temporarily consider the derived Satake isomorphism
    \[
    \Sph_{G,x} \simeq
    \Sph_{\LG,x}^\spec \simeq \ICoh_\Nilp\left( ( \pt \underset{\Lg}{\times} \pt ) /\LG \right).
    \]
The na\"ive diagram
\[
\begin{tikzcd}[
  column sep={8em,between origins},
  row sep={4em,between origins}
]
    \Sph_{G,x} \arrow[r, "\Sat"] \arrow[d, "\sigma" left] \arrow[dr, phantom, "(\tau)"] & \Sph_{\LG,x}^\spec \arrow[d, "\sigma"] \\
    \Sph_{G,x}^\rev \arrow[r, "\Sat^\rev" below] & \Sph_{\LG,x}^{\spec, \rev}
\end{tikzcd}
\]
of algbera isomorphisms is not commutative: the commutativity is only restored after inserting the Cartan involution. Note that the diagram
\begin{equation}
\label{eqn:Sat_Cartan}
\begin{tikzcd}[
  column sep={7em,between origins},
  row sep={4em,between origins}
]
    \Rep(\LG) \arrow[r] \arrow[d, "\id" left] & \Sph_{\LG,x}^\spec \arrow[d, "\sigma"] \\
    \Rep(\LG)^\rev \arrow[r,] & \Sph_{\LG,x}^{\spec, \rev}
\end{tikzcd}
\end{equation}
is canonically commutative because the horizonal map is given by pushforward along the diagonal
    \[\BB \LG \to \BB\LG \underset{\Lg/\LG}{\times} \BB\LG.\]
Since Serre duality on $\ICoh_\Nilp(\LS_\LG)$ is naturally $\Rep(\LG)$-linear, the Cartan involution in (\ref{eqn:Sat_Cartan}) explains the presence of the Cartan involution in (\ref{eqn:Serre_PsId_no_shift}). Similarly, the Cartan involution also appears in the Whittaker duality:
\[
\begin{tikzcd}[
  column sep={8em,between origins},
  row sep={4em,between origins}
]
  \Whit_\psi(G)_x \arrow[d, "\DD{[-2\delta_N]}" left] \arrow[dr, phantom, "\tau"] & \Rep(\LG)_x \arrow[d] \arrow[l, "\CS" above]\\
  \Whit_{-\psi}(G)_x^\vee \arrow[r, "\CS^\vee"] & \Rep(\LG)_x^\vee
\end{tikzcd}
\]

\subsubsection{} Applying the expected compatibilities to every object appearing Theorem \ref{thm:Ran}, we obtain the diagram
\[
\begin{tikzcd}[
  column sep={10em,between origins},
  row sep={4em,between origins}
]
      \Rep(\LG)_x^\vee \arrow[d, "\ev_{x,*}^\vee" left] & \Rep(\LG)_x \arrow[l] \arrow[d, "\ev_x^*"] \\
      \ICoh_\Nilp(\LS_\LG)^\vee & \ICoh_\Nilp(\LS_\LG) \arrow[l, "\DD{[2\delta_G]}"]
\end{tikzcd}
\]

\subsection{Outline of the proof}

\subsubsection{Notation} Let us write $(-)_\co$ for the Verdier duality $\DD(-)$ and
    \[\PsId_\Whit : \Whit_{-\psi}(\Bun_N)^\vee \to \Whit_\psi(\Bun_N)\]
for its inverse.

\subsubsection{} Our first task will be the construction of a natural transformation
    \[\eta_G : \Poinc_! \cdot \PsId_\Whit \to \PsId_! \cdot \coeff_*^\vee \, [2\delta_G].\]
This problem is actually very concrete: the data of $\eta_G$ is equivalent to a morphism
    \[
    \Poinc_! \cdot \WW_\vac \to \PsId_! \cdot \coeff_*^\vee \cdot \WW_{\vac, \co} \, [2 \delta_G]
    \]
between two objects of $\DMod(\Bun_G)$. We will perform this construction in \S\ref{section:construction}. Once this has been done, we will have to prove that $\eta_G$ is an isomorphism.

\subsubsection{} Our proof that $\eta_G$ is an isomorphism is based on the decomposition of $\DMod(\Bun_G)$ into its cuspidal and Eisenstein parts. Recall that the full subcategory $\DMod(\Bun_G)_\cusp \mon \DMod(\Bun_G)$ of cuspidal objects has the following equivalent characterizations:
\begin{itemize}
    \item It consists of objects $\FF$ which are annihilated by $\CT_{P,*}$ for every proper parabolic $P$.
    \item It consists of objects $\FF$ which are annihilated by $\CT_{P,!}$ for every proper parabolic $P$.
    \item It is the right orthogonal to $\DMod(\Bun_G)_\Eis$.
\end{itemize}
Here the full subcategory $\DMod(\Bun_G)_\Eis \mon \DMod(\Bun_G)$ is characterized by the following equivalent characterizations:
\begin{itemize}
    \item It is the full subcategory generated by the images of $\Eis_{P,!}$ for all proper parabolics $P$.
    \item It is the left orthogonal to $\DMod(\Bun_G)_\cusp$.
\end{itemize}
The category $\DMod(\Bun_G)_\Eis$ contains but is not usually generated by the images of $\Eis_{P,*}$ for all proper parabolics $P$. We will prove Theorem \ref{thm:vac} by establishing the following:

\begin{lemma}
\label{lemma:cofib_eta_Eis}
The cofiber of $\eta_G$ is Eisenstein series.
\end{lemma}

\begin{lemma}
\label{lemma:cofib_eta_cusp}
The cofiber of $\eta_G$ is cuspidal.
\end{lemma}

\subsection{Comparison with the na\"ive pseudo-identity}

\subsubsection{} Recall that there is a na\"ive variant of the pseudo-identity functor, denoted
    \[\PsId_* : \DMod(\Bun_G)^\vee \to \DMod(\Bun_G),\]
corresponding to the object
    \[\Delta_* \cdot \omega_{\Bun_G} \in \DMod(\Bun_G \times \Bun_G).\]
It is related to $\PsId_!$ by a canonical map
    \[
    \PsId_! \, [2 \delta_G + z_G] \to \PsId_*
    \qquad
    z_G = \dim Z_G
    \]
whose cofiber is Eisenstein series.

\subsubsection{} The functor $\Poinc_!$ has a variant
    \[\Poinc_* = \pi_* \cdot \oblv,\]
related to $\Poinc_!$ by a canonical map
\begin{equation}
    \label{eqn:Poinc*!2}
    \Poinc_! \, [z_G] \to \Poinc_*.
\end{equation}
The version of our main theorem with $\PsId_*$ and $\Poinc_*$ is actually trivial:

\begin{proposition}
\label{proposition:trivial_thm}
There is a canonical commutative diagram
\[
\begin{tikzcd}[
  column sep={10em,between origins},
  row sep={4em,between origins}
]
    \Whit_{-\psi}(\Bun_N)^\vee \arrow[d, "\coeff_*^\vee" left] \arrow[r, "\PsId"] & \Whit_\psi(\Bun_N) \arrow[d, "\Poinc_*"] \\
    \DMod(\Bun_G)^\vee \arrow[r, "\PsId_*" below] & \DMod(\Bun_G)
\end{tikzcd}
\]
\end{proposition}

\subsubsection{Proof} Performing base change on the diagram
\[
\begin{tikzcd}[
  column sep={4em,between origins},
  row sep={4em,between origins}
]
    & \Bun_N \arrow[dl] \arrow[dr] & & \\
    \Bun_G \arrow[dr, "\Delta"] & & \Bun_N \times \Bun_G \arrow[dl, "\pi \times \id"] \arrow[dr, "\pi_2"] & \\
    & \Bun_G \times \Bun_G & & \Bun_G
\end{tikzcd}
\]
gives
\[
    \PsId_* \cdot \coeff_*^\vee \cdot \WW_\co
    \simeq \pi_{2,*} \cdot \left\{ \WW \shimes (\pi \times \id)^! \cdot \Delta_* \cdot \omega_{\Bun_G} \right\}
    \simeq \Poinc_* \cdot \WW.
\]

\subsubsection{} Combining $\eta_G$ with the morphism relating $\PsId_!$ to $\PsId_*$, we obtain a map
\begin{equation}
\label{eqn:Poinc*!1}
\begin{split}
  \Poinc_! \, [z_G]
  & \to \PsId_! \cdot \coeff^\vee \cdot \PsId_\Whit^{-1} \, [2 \delta_G + z_G] \\
  & \qquad \qquad \to \PsId_* \cdot \coeff^\vee \cdot \PsId_\Whit^{-1} \simeq \Poinc_*.
\end{split}
\end{equation}
Then Lemma \ref{lemma:cofib_eta_Eis} is equivalent to the assertion that the cofiber of this map is Eisenstein series. We will prove this in \S\ref{section:cofib_Eis} by verifying the following two assertions:

\begin{proposition}
\label{proposition:poinc_identification}
The map (\ref{eqn:Poinc*!1}) canonically identifies with (\ref{eqn:Poinc*!2}).
\end{proposition}

\begin{proposition}
\label{proposition:poinc_principal_series}
The cofiber of (\ref{eqn:Poinc*!2}) is Eisenstein series.
\end{proposition}

\subsubsection{Remark} Combining Theorem \ref{thm:vac} with Proposition \ref{proposition:trivial_thm} gives
    \[
    \Poinc_* \simeq \DL \cdot \Poinc_!
    \]
where
    \[\DL = \PsId_* \cdot \PsId_!^{-1} \, [-2\delta_G]\]
is the Deligne--Lusztig functor. This means that $\Poinc_* \cdot \WW_\vac$ should correspond under Langlands duality to the Steinberg object constructed by D. Beraldo \cite{dario_DL}.

\subsection{Jacquet functors} 

\subsubsection{} In order to show that the cofiber of $\eta_G$ is cuspidal, we need to prove that the morphism
    \[
    \CT_{P, !}^- \cdot \Poinc_! \cdot \WW_\vac 
    \to \CT_{P, !}^- \cdot \PsId_! \cdot \coeff^\vee \cdot \WW_{\vac,\co} \, [2\delta_G]
    \]
is an equivalence for every proper parabolic $P$. We will do this by writing both sides in terms of Poincare series for $M$ and applying induction. This will force us to consider the Ran version of the Whittaker category.

\subsubsection{The LHS} We will construct a Jacquet functor 
  \[\Jac : \Whit(\Bun_N) \to \Whit(M)_\Ran\]
fitting into a commutative square:
\[
\begin{tikzcd}[
  column sep={7.5em,between origins},
  row sep={4em,between origins}
]
    \Whit(G) \arrow[d, "\Poinc_!" left] \arrow[r, "\Jac"] & \Whit(M)_\Ran \arrow[d, "\Poinc_!"] \\
    \DMod(\Bun_G) \arrow[r, "\CT_P^-"] & \DMod(\Bun_M) 
\end{tikzcd}
\]
The right vertical functor appearing in this diagram is the composite
    \[\Whit(M)_\Ran \to \DMod(\Bun_M)_\Ran \to \DMod(\Bun_M),\]
where the first arrow is the family version of Poincare series appearing in Theorem \ref{thm:Ran} for $M$, and the second arrow is integration over the Ran space. The proof that this square is commutative will require geometric input in the form of a ULA property established in \S\ref{section:ULA}.

\subsubsection{The RHS} First recall that
    \[
    \CT^-_{P, !} \cdot \PsId_{G,!} 
    \simeq \PsId_{M,!} \cdot (\Eis_{P,!}^-)^\vee
    \]
because, as objects of $\DMod(\Bun_G \times \Bun_M)$, both identify by base-change with
    \[(p \times q)_! \cdot k_{\Bun_{P^-}}.\]
Then, similarly to the above, we will construct a functor
    \[\Jac_\co : \Whit(M)_\Ran \to \Whit(\Bun_N)\]
equipped with a commutative square:
\[
\begin{tikzcd}[
  column sep={8em,between origins},
  row sep={4em,between origins}
]
  \DMod(\Bun_M) \arrow[r, "\Eis_!^-"] \arrow[d, "\coeff_*" left] & \DMod(\Bun_G) \arrow[d, "\coeff_*"] \\
  \Whit(M)_\Ran \arrow[r, "\Jac_{\co}"] & \Whit(G)
\end{tikzcd}
\]
The functors $\Jac$ and $\Jac_\co$ will be dual to each other in the sense that there a commutative diagram:
\[
\begin{tikzcd}[
  column sep={8em,between origins},
  row sep={4em,between origins}
]
  \Whit(G) \arrow[r, "\Jac"] \arrow[d, "\DD{[2\delta_G]}" left] & \Whit(M)_\Ran \arrow[d, "\DD{[2\delta_M]}"] \\
  \Whit(G)^\vee \arrow[r, "\Jac_{\co}^\vee"] & \Whit(M)_\Ran^\vee 
\end{tikzcd}
\]

\subsubsection{} We organize this discussion by pasting everything into a large diagram.
\begin{equation}
\label{eqn:the_cube}
\begin{tikzcd}[
  column sep={5em,between origins},
  row sep={4em,between origins}
]
  \Whit(M)_\Ran^\vee \arrow[dddd, "\coeff_*^\vee" left] & & & & \Whit(M)_\Ran \arrow[llll, "\DD{[2\delta_M]}" above] \arrow[dddd, "\Poinc_!"] \\
  & \Whit(G)^\vee \arrow[dd, "\coeff_*^\vee" left] \arrow[ul, "\Jac_\co^\vee" above right] & & \Whit(G) \arrow[ll, "\DD{[2\delta_G]}" above] \arrow[dd, "\Poinc_!"] \arrow[ur, "\Jac"] & \\
  & & & & \\
  & \DMod(\Bun_G)^\vee \arrow[dl, "\Eis_{P,!}^{-,\vee}"] \arrow[rr, "\PsId_!"] & & \DMod(\Bun_G) \arrow[dr, "\CT_{P,!}^-" below left]  &  \\
  \DMod(\Bun_M)^\vee \arrow[rrrr, "\PsId_{M,!}" below] & & & & \DMod(\Bun_M)
\end{tikzcd}
\end{equation}
In \S\ref{section:Ran_upgrades}, we shall extend the construction of $\eta_G$ to a $\Sph_{G,\Ran}$-linear morphism
    \[
    \eta_{G,\Ran} : \Poinc_{\Ran, !} \cdot \PsId_{\Ran, !}^\Whit \to \PsId_{G,!} \cdot \coeff_{\Ran, *}^\vee \, [2\delta_G].
    \]
It is compatible in the sense that its evaluation on $\WW_\vac$ agrees with $\eta_G$. Therefore, Theorem \ref{thm:vac} for a given group implies Theorem \ref{thm:Ran} for that same group.

\begin{lemma}
\label{lemma:the_cube_commutes}
Complete (\ref{eqn:the_cube}) into a hexahedron of functors and natural transformations by pasting $\eta_{M,\Ran}$ along the outer edges. Then there exists a 3-morphism rendering the cube commutative.
\end{lemma}

\subsubsection{} We were not able to prove Lemma \ref{lemma:the_cube_commutes} directly. Instead, the proof of the theorem goes through a slightly different statement from which this lemma can be deduced. This is a rather unsatisfactory situation because the resulting 3-morphism is non-canonical.

\subsubsection{} Let us explain how Theorem \ref{thm:vac} would follow by induction on $G$ if we knew Lemma \ref{lemma:the_cube_commutes}. Suppose that $\eta_M$ is an isomorphism for all proper Levi subgroups of $G$. By $\Sph_{M,\Ran}$-linearity, $\eta_{M,\Ran}$ is a an isomorphism for all proper Levis $M$. Then Lemma \ref{lemma:the_cube_commutes} implies that $\CT_{P,!}^- \cdot \eta_G$ is an isomorphism for all proper Levis $M$. This proves Lemma \ref{lemma:cofib_eta_cusp} for $G$. Combining this with Lemma \ref{lemma:cofib_eta_Eis} proves that $\eta_G$ is an isomorphism.

\subsection{Acknowledgements}

\subsubsection{} This problem was suggested by our advisor D. Gaitsgory, and this work would not have existed without his steady guidance and support. We are grateful to S. Raskin for conversations which lead to the argument of \S\ref{section:ULA}. We are also grateful to L. Chen and Y. Fu for many patient explanations and for their encouragement.

\subsubsection{} We thank L. Chen and D. Gaitsgory for their comments on an earlier draft of this work.

%% file: construction.tex
\section{Construction of the map}
\label{section:construction}

\subsubsection{} In this section we will construct the map
  \[\eta_G : \Poinc_! \cdot \PsId_\Whit^{-1} \to \PsId_! \cdot \coeff_*^\vee [2 \delta_G].\]

\subsection{Reminders on Vinberg's semigroup}

\subsubsection{} We write
    \[T_\ad = T/Z_G \simeq (\Gm)^\Delta\]
for the adjoint torus. We write
    \[T_\ad^+ = (\AA^1)^\Delta,\]
viewed as a semigroup. The points in $T_\ad^+$ with zero-one co\"ordinates are in natural bijection with subsets of $\Delta$ and hence standard parabolics. We write $c_P$ for the point corresponding to the parabolic $P$. The $T$ orbit of $c_P$ identifies with $(T_M)_\ad \simeq T/Z_M$, and these orbits form a stratification of $T_\ad^+$.

\subsubsection{Vinberg's semigroup} Consider the group homomorphism
    \[G \overset{Z_G}{\times} T \to T_\ad.\]
This map is naturally $G \times T \times G$ equivariant, where $G \times G$ acts on the left and right on the source and trivially on the target. Vinberg's semigroup $\Vin_G$ is an algebraic semigroup equipped with a $G \times G \times T$-equivariant semigroup morphism
    \[\Vin_G \to T_\ad^+\]
extending the above group homomorphism. This morphism has a canonical ``unit section''
    \[s : T^+_\ad \to \Vin_G\]
whose orbit under $G \times G$ is an open subscheme denoted $\Vin_G^\circ$.

\subsubsection{Stratification} The fiber of
    \[T_\ad^+ \to \Vin_G^\circ \mon \Vin_G\]
over $c_P$ is canonically identified with
\begin{equation}
\label{eqn:Vin_G_strata}    
    \pt \simeq P^-/U_P^- \overset{M}{\times} U_P \bs P
    \to G/U_P^- \overset{M}{\times} U_P \bs G 
    \mon \overline{G/U_P^- \overset{M}{\times} U_P \bs G}
\end{equation}
with the natural $G \times G$ action. The residual action of $Z_M$ is descended from the right and left actions on $G/U_P^-$ and $U_P \bs G$, respectively.

\subsubsection{Generic maps} Let us recall the notion of generic maps from [J. Barlev]. For any scheme $S$, an open subset $U \mon X \times S$ is called a domain if it is open dense after any base change $S' \to S$. For any target $\YY$, we define the prestack
    \[
    \Maps^\gen(X, \YY) : S \mapsto
    \lim_U \Maps(U, \YY)
    \]
where the limit is taken over all domains $U \mon X \times S$. Given a map $\YY_1 \to \YY_2$, we define $\Maps^\gen(X, \YY_1 \to \YY_2)$ as the fiber product:
\[
\begin{tikzcd}
    \Maps^\gen(X, \YY_1 \to \YY_2) \arrow[r] \arrow[d] &  \Maps(X, \YY_2) \arrow[d] \\
    \Maps^\gen(X, \YY_1) \arrow[r] & \Maps^\gen(X, \YY_1)
\end{tikzcd}
\]

\subsubsection{} We define
    \[
    \Vin\Bun_G = \Maps^\gen \left(X, G \bs \Vin^\circ_G / G \mon G \bs \Vin_G / G \right).
    \]
This a stack over $\Bun_G \times \Bun_G \times T_\ad^+$, and we may take the quotient by $T$ to obtain
    \[\Bunc_G = \VBun_G / T,\]
which is a stack over $\Bun_G \times \Bun_G \times T_\ad^+/T$. The stratification of $T_\ad^+/T$ by the poset of standard parabolics induces a stratification of $\Bunc_G$, whose pieces we shall denote $\Bunc_{G,P}$. Note that the open stratum identifies with
    \[\Bunc_{G,G} \simeq \Bun_G \times \BB Z_G.\]

\subsubsection{Reminders on $\PsId_*$ and $\PsId_!$} There is a sequence of maps
  \[\Bun_G \overset{r}{\longrightarrow} \Bun_G \times \BB Z_G \overset{j}{\hookrightarrow} \Bunc_G \overset{\Deltac}{\longrightarrow} \Bun_G \times \Bun_G\]
with the following features:
\begin{itemize}
    \item We have $r_! \simeq r_* [z_G]$.
    \item The map $\Deltac$ is proper.
\end{itemize}
Therefore, we obtain a map
    \[\Delta_! \cdot k_{\Bun_G} \to \Delta_* \cdot k_{\Bun_G} [z_G] \simeq \Delta_* \cdot \omega_{\Bun_G} [z_G - 2 \delta_G].\]
This gives the map
    \[\PsId_! \to \PsId_* [z_G - 2 \delta_G].\]

\subsection{Definition of $\eta_G$}

\subsubsection{Notation} We will write
    \[r_{G,!} = j_! \cdot r_! \cdot k_{\Bun_G} \in \DMod\left( \Bunc_G \right).\]

\subsubsection{} The object $\WW_\vac$ extends cleanly along $\Bun_N \mon \Bunc_N$. We use the same symbol for this extension, and we make the identification
    \[\Whit\left(\Bunc_N\right) \simeq \Whit(\Bun_N).\]
We shall write the definition of $\eta_G$ using $\Bunc_N$ even though it is slightly more complicated because the $\Bunc_N$ version will be needed in the construction of the Jacquet functor.

\subsubsection{} Consider the Cartesian diagram below.
\[
\begin{tikzcd}[
  column sep={4em,between origins},
  row sep={4em,between origins}
]
    & \Bunc_N \arrow[dl] \arrow[dr] & & & \\    
  \Bun_G \arrow[dr] & & \Bunc_N \underset{\Bun_G}{\times} \Bunc_G \arrow[dl] \arrow[dr, "\aph\Deltac"] & & \\
    & \Bunc_G \arrow[dr] & & \Bunc_N \times \Bun_G \arrow[dl] \arrow[dr, "\pr_2"] & \\
    &  & \Bun_G \times \Bun_G & & \Bun_G
\end{tikzcd}
\]
We write
    \[\pi_2 = \pr_2 \cdot \aph\Deltac.\]
Note that we have the identifications
    \[
    \Poinc_! \cdot \WW_\vac \simeq \pi_{2,!} \left( \WW_\vac \stoximes r_{G,!} \right)
    \qquad
    \PsId_! \cdot \coeff_*^\vee \cdot \WW_{\vac,\co} \simeq \pi_{2,*} \left( \WW_\vac \shoximes r_{G,!} \right)
    \]
by base change. Therefore, we expect that the purity morphism and forgetting supports combine to give a map
    \[
    \pi_{2,!} \left( \WW_\vac \shoximes r_{G,!} \right)
    \to \pi_{2,*} \left( \WW_\vac \stoximes r_{G,!} \right) [2 \delta_G].
    \]
Unfortunately, $\pi_2$ is not representable, so there is no forgetting supports morphism $\pi_{2,!} \to \pi_{2,*}$. To circumvent this issue, we will construct an open substack
\begin{equation}
\label{eqn:trans_open}
  j : \left( \Bunc_N \!\underset{\Bun_G}{\times}\! \Bunc_G \right)^\trans \mon \Bunc_N \!\underset{\Bun_G}{\times}\! \Bunc_G
\end{equation}
containing $\Bunc_N \times \BB Z_G$ with the following features:

\begin{proposition}
\label{proposition:cis_vanishing}
Let
\[
   i: \left( \Bunc_N \!\underset{\Bun_G}{\times}\! \Bunc_G \right)^\cis \mon \Bunc_N \!\underset{\Bun_G}{\times}\! \Bunc_G
\]
denote the closed complement to $j$. Then
  \[\pi_{2,*} \cdot i_* \cdot i^! \left( \WW \shoximes r_{G,!} \right) \simeq 0.\]
\end{proposition}

\begin{proposition}
\label{proposition:trans_rep}
The projection $\pi_2 \cdot j$ onto $\Bun_G$ is representable.
\end{proposition}

\subsubsection{} We think of this open substack a locus of transversality, and we will write $(-)^\trans$ for restriction along $j$. Then $\eta_G$ is defined via the diagram:
\[
\begin{tikzcd}[
  column sep={6em,between origins},
  row sep={4em,between origins}
]
  & \pi_{2, !} \cdot \left( \WW \stoximes r_! \right)^\trans \arrow[dr] \arrow[dl, "\sim"] & & \pi_{2, *} \cdot \left( \WW \shoximes r_! \right) [2 \delta_G] \arrow[dl, "\sim"] \\
  \pi_{2, !} \cdot \left( \WW \stoximes r_! \right) & & \pi_{2, *} \cdot \left( \WW \shoximes r_! \right)^\trans [2 \delta_G] &
\end{tikzcd}
\]
The left diagonal map is an isomorphism because $\WW \stoximes r_!$ is $!$-extended from $\Bun_N \times \BB Z_G$. The right diagonal map is an isomorphism by virtue of Proposition \ref{proposition:cis_vanishing}. The middle diagonal map is a combination of the purity morphism and the forgetting supports morphism for $\pi_2 \cdot j$, which makes sense by Proposition \ref{proposition:trans_rep}.

\subsubsection{} We now give the definition of (\ref{eqn:trans_open}). The proofs of Propositions \ref{proposition:cis_vanishing} and \ref{proposition:trans_rep} will appear in the subsequent subsections.

\begin{proposition}
\label{proposition:open_orbit}
The $N \times G$ orbit of $s$ is an open subscheme $\Vin_G^{\circ\circ} \mon \Vin_G^\circ$. It is equal to the $B \times G$ orbit of $s$.
\end{proposition}

\subsubsection{Proof of Proposition \ref{proposition:open_orbit}} Since smooth morphisms are open, it suffices to show that the morphism
  \[N \times T_\ad^+ \times G \to N \times \Vin_G^\circ \times G \to \Vin_G^\circ\]
given by acting on the unit section is smooth. This is a map of schemes smooth over $T_\ad^+$, so we can check the smoothness fiberwise using the explicit description (\ref{eqn:Vin_G_strata}) of the fibers.

\subsubsection{} Now
    \[
    \Bunc_B \underset{\Bun_G}{\times} \Bunc_G
    \simeq \Maps^\gen\left(X, B \bs \Vin_G^\circ / G \mon T \bs \overline{N\bs G} \overset{G}{\times} \Vin_G / G \right)
    \]
and we define
    \[
    \left( \Bunc_B \underset{\Bun_G}{\times} \Bunc_G \right)^\trans
    = \Maps^\gen\left(X, B \bs \Vin_G^{\circ\circ} / G \mon T \bs \overline{N\bs G} \overset{G}{\times} \Vin_G / G \right)
    \]
We take
    \[
    \left( \Bunc_N \underset{\Bun_G}{\times} \Bunc_G \right)^\trans
    = \Bunc_N \underset{\Bunc_B}{\times} \left( \Bunc_B \underset{\Bun_G}{\times} \Bunc_G \right)^\trans
    \]
The map
    \[
    \Bunc_N \times \BB Z_G \to \Bunc_N \underset{\Bun_G}{\times} \Bunc_G
    \]
factors through the transverse locus because $\Vin_G^{\circ\circ}$ contains $\Vin_G^\circ \underset{T^+_\ad}{\times} T_\ad$.

\subsubsection{Sanity check} Let us study the extent to which this device depends on the choice of open substack. Suppose that we have a diagram
\[
\begin{tikzcd}
  \UU \arrow[r, hook, "j"] & \YY \arrow[d, "\pi"] \\
  & S
\end{tikzcd}
\]
in which $j$ is an open embedding and $\pi$ is representable. Then for any map $\FF \to \GG$ of sheaves on $\YY$, we have an obvious commutative diagram:
\[
\begin{tikzcd}[
  column sep={5em,between origins},
  row sep={4em,between origins}
]
  j_! \cdot \FF_\UU \arrow[r] \arrow[d] & \FF \arrow[r] \arrow[d] & j_* \cdot \FF_\UU \arrow[d] \\
  j_! \cdot \GG_\UU \arrow[r] & \GG \arrow[r] & j_* \cdot \GG_\UU
\end{tikzcd}
\]
This induces
\[
\begin{tikzcd}[
  column sep={5em,between origins},
  row sep={1.5em,between origins}
]
  \pi_! \cdot j_! \cdot \FF_\UU \arrow[r] \arrow[dd] & \pi_! \cdot \FF \arrow[dr] \arrow[ddd] & \\
   & & \pi_* \cdot j_* \cdot \FF_\UU \arrow[dd] \\
  \pi_! \cdot j_! \cdot \GG_\UU \arrow[dr] & & \\
   & \pi_* \cdot \GG \arrow[r] & \pi_* \cdot j_* \cdot \GG_\UU
\end{tikzcd}
\]
This diagram exhibits the sense in which our morphism does not depend on the choice of $\UU$. 

\subsubsection{} This sanity check shows that the definition of $\eta_G$ would be unchanged if we replaced $\Bunc_N$ with $\Bun_N$ everywhere.

\subsection{Example: what happens for $\SL_2$}

\subsubsection{} The rest of this paper has no logical dependencies on this subsection.

\subsubsection{} We fix $G = \SL_2$ for this subsection. In this case the objects at hand admit the following descriptions:
\begin{itemize}
    \item The stack $\Vin\Bun_G$ classifies $(\EE_1, \EE_2, \alpha)$ where the $\EE_i$ are a pair of $\SL_2$-bundles and $\alpha$ is a nonzero map of coherent sheaves
        \[\alpha : \EE_1 \to \EE_2.\]
    The torus $T \simeq \Gm$ acts by scaling $\alpha$.
    \item The quotient $\Bunc_G$ classifies $(\EE_1, \EE_2, L, \alpha)$ where the $\EE_i$ are a pair of $\SL_2$-bundles, $L$ is a line, and $\alpha$ is a nonzero map of coherent sheaves
        \[\alpha : \EE_1 \to \EE_2 \otimes L.\]
    The open stratum $\Bunc_{G,G}$ is the locus where $\alpha$ is an isomorphism, and the closed stratum $\Bunc_{G,B}$ is the locus where $\alpha$ has rank one.
    \item The fiber product $\Bunc_N \underset{\Bun_G}{\times} \Bunc_G$ classifies the data of
        \[\omega^{1/2} \mon \EE_1 \to \EE_2 \otimes L.\]
    The transverse locus is the open substack where the composed map $\omega^{1/2} \to \EE_2 \otimes L$ is nonzero.
\end{itemize}

\subsubsection{} Let's verify by hand that the transverse locus is representable over $\Bun_G$. The $G$-stratum of the transverse locus is simply
    \[\Bunc_N \times \BB Z_G \to \Bun_G\]
so we may restrict attention to the $B$-stratum. This stratum has a further stratification by the degree of $\im(\alpha)$. Now the stratum where $\im(\alpha)$ has degree $d$ can be perceived as classifying a point $\EE_2$ of $\Bun_G$ along with:
\begin{itemize}
    \item[(i)] a choice of degree $d$ divisor $D$,
    \item[(ii)] a line $L$ along with an injective map of coherent sheaves
        \[\omega^{1/2}(D) \mon \EE_2 \otimes L,\]
    and
    \item[(iii)] an $\SL_2$ bundle $\EE_1$ fitting into a triangle 
        \[
        \begin{tikzcd}[
          column sep={1em},
          row sep={1em}
        ]
            \omega^{1/2} \arrow[d, hook] \arrow[dr, hook] & \\
            \EE_1 \arrow[r, ->>] & \omega^{1/2}(D)
        \end{tikzcd}
        \]
\end{itemize}
We claim that the stacks parameterizing each additional data are schematic over the previous ones. So for (i) we are looking at the projection
    \[\Sym^d(X) \times \Bun_G \to \Bun_G,\]
which is schematic because $\Sym^d(X)$ is a scheme. Then the choices in (ii) are parameterized by the generalized projective bundle over $\Sym^d(X) \times \Bun_G$ whose fibers are 
    \[\bfP H^0\left( \Maps(\omega^{1/2}(D), \EE_1) \right).\]
Finally, the data in (iii) are parametrized by the vector bundle over this projective bundle whose fibers are
\begin{align*}
    & \fib \left\{ \Maps\left(\omega^{1/2}(D), \omega^{-1/2}(-D)[1] \right) \to \Maps\left( \omega^{1/2}, \omega^{-1/2}(-D)[1] \right) \right\}   \\
    & \qquad \qquad \simeq \Maps\left(\omega^{1/2}_D, \omega^{-1/2}(-D) [1] \right)
    \simeq \Maps\left( \OO, (\omega_D)^{\otimes 2} \right)^\vee.
\end{align*}

\subsection{Proof of Proposition \ref{proposition:cis_vanishing}}

\subsubsection{} We prove Proposition \ref{proposition:cis_vanishing} by a series of reductions. Let us write
    \[\FF = \WW \shoximes r_{G,!}.\]

\subsubsection{First reduction} We use the stratification of $\Bunc_G$ into the strata $\Bunc_{G,P}$. It suffices to show that the !-restrictions of $\FF$ to each
  \[\left(\Bunc_N \underset{\Bun_G}{\times} \Bunc_{G,P}\right)^\cis\]
are killed by $*$-pushforward to $\Bun_G$.

\subsubsection{A generality} Let $T$ be any algebraic torus, and consider a principal $T$-bundle
\[
\begin{tikzcd}[
  column sep={4em,between origins},
  row sep={3em,between origins}
]
    P \arrow[d, "f" left] \arrow[r] & \pt \arrow[d] \\
    S \arrow[r, "\alpha"] & \BB T
\end{tikzcd}
\]
For any sheaf $\FF$ on $S$ we have
    \[\Gamma(S, \FF) \simeq \Maps(k_{\BB T}, \alpha_* \cdot \FF) \simeq \Maps_{\rH_*(T)}\left(r^! \cdot k_{\BB T}, \Gamma(P, \FF) \right).\]

\subsubsection{Second reduction} Let $\Vin\Bun_{G,P}$ denote the fiber of $\Vin\Bun_G$ over the point $c_P \in T_\ad^+$. We have
    \[\Bun_{G,P} \simeq \Vin\Bun_{G,P}/Z_M.\]
By the above, it suffices to show that the !-pullbacks of $\FF$ to
  \[\left( \Bunc_N \underset{\Bun_G}{\times} \VBun_{G,P} \right)^\cis\]
are killed by $*$-pushforward to $\Bun_G$.

\subsubsection{Third reduction} Recall that $\VBun_{G,P}$ has a stratification of the form
\begin{equation}
\label{eqn:VinBun_G_str}
    \Bun_{P^-} \underset{\Bun_M}{\times} \HH_M^+ \underset{\Bun_M}{\times} \Bun_P \to \VBun_{G,P}.      
\end{equation}
Therefore, it suffices to show that the !-pullback of $\FF$ to this stratification is killed by $*$-pushforward along
  \[\left( \Bunc_N \underset{\Bun_G}{\times} \Bun_{P^-} \underset{\Bun_M}{\times} \HH_M^+ \underset{\Bun_M}{\times} \Bun_P \right)^\cis \to \Bun_P \to \Bun_G.\]

\subsubsection{Fourth reduction} Finally, we recall that  
  \[\Bunc_N \underset{\Bun_G}{\times} \Bun_{P^-}\]
has a relative position stratification by
  \[W/W_P.\]
It suffices to prove the analogous statement for
  \[
  \pi_2^\gr : \left( \Bunc_N \underset{\Bun_G}{\times} \Bun_{P^-} \right)^w \underset{\Bun_M}{\times} \HH_M^+ \underset{\Bun_M}{\times} \Bun_P  \to \HH_M^+ \underset{\Bun_M}{\times} \Bun_P
  \]
where $w \neq 1$.

\subsubsection{End of the proof} By \cite[Lemma 4.2.1]{lin_nearby_cycles}, the !-restriction of $r_!$ to the above stack is !-pulled back along $\pi_2^\gr$. Therefore, it suffices to observe that $\WW$ is annihilated by $*$-pushforward along
    \[
    \left( \Bunc_N \underset{\Bun_G}{\times} \Bun_{P^-} \right)^w 
    \to \Bun_M.
    \]
This is a standard fact whose proof we review in \S\ref{subsection:Jacquet_Whittaker}.

\subsection{Proof of Proposition \ref{proposition:trans_rep}}

\subsubsection{} To simplify notation, we will prove the analogue of Proposition \ref{proposition:trans_rep} for the \emph{untwisted} version of $\Bun_N$. The twisted version can be proved in the same way by using the appropriate non-constant group scheme twist of $N$ and replacing the notion of generic maps with generic sections.

\subsubsection{Proof of Proposition \ref{proposition:trans_rep}} Recall that an Artin stack is an algebraic space if the stabilizer group of every geometric point is trivial. Therefore, we only need to check that our map induces injections on stabilizer groups when evaluated at geometric points (note that this property is stable under base change). We shall do this after passage to a stratification

\subsubsection{Stratification} We will again use the stratification of $\Bunc_G$ into $\Bunc_{G,P}$. We further stratify $\Bunc_{G,P}$ using the quotient of (\ref{eqn:VinBun_G_str}) by $Z_M$. Therefore, it suffices to show that the map
    \[
    \left(\Bunc_N \underset{\Bun_G}{\times} \Bun_{P^-} \right)^\trans \underset{\Bun_M}{\times} \HH_M^+/Z_M \underset{\Bun_M}{\times} \Bun_P \to \Bun_G
    \]
is representable. Since the morphism $\Bun_P \to \Bun_G$ is schematic, we only need to check that the map
    \[
    \alpha : \left(\Bunc_N \underset{\Bun_G}{\times} \Bun_{P^-} \right)^\trans \underset{\Bun_M}{\times} \HH_M^+/Z_M \to \Bun_M
    \]
is representable. 

\subsubsection{} Now recall that the source of the stratification (\ref{eqn:VinBun_G_str}) is given by
    \[
    \Maps^\gen\left(X, P^- \bs M/P \mon P^- \bs M^+ /P \right)
    \]
for a certain monoid $M^+$ containing $M$ as its unit group. Therefore, in terms of generic maps, the morphism $\alpha$ can be written as the composite:
\begin{align*}
    & \Maps^\gen\left( X, \BB N_M \mon N \bs G \overset{P^-}{\times} M^+/M \right)/Z_M \\
    & \qquad \mon \Maps^\gen\left( X, \BB (N_M \cdot Z_M) \mon N \bs G \overset{P^-}{\times} M^+/M \cdot Z_M \right) \\
    & \qquad \qquad \overset{\beta}{\longrightarrow} \Maps(X, \BB M)
\end{align*}
The first arrow is a closed immersion because it is a base change of $\BB M \to \Bun_M$. Therefore, it suffices that to show that $\beta$ induces injections on stabilizer groups when evaluated at geometric points. 

\subsubsection{} Let us write
    \[
    \YY = N \bs G \overset{P^-}{\times} M^+/Z_M \cdot M
    \qquad
    V = N_M \cdot Z_M
    \]
for convenience. Then $\beta$ fits into the following Cartesian diagram. 
\[
\begin{tikzcd}[
  column sep={12em,between origins},
  row sep={4em,between origins}
]
    \Maps^\gen\left(X, \BB V \mon \YY \right) \arrow[r] \arrow[d] & \Maps(X, \YY) \arrow[d, "\star"] & \\
    \Maps^\gen\left(X, \BB V \to \BB M \right) \arrow[r] \arrow[d] & \Maps^\gen\left(X, \YY \to \BB M \right) \arrow[r] \arrow[d] & \Maps(X, \BB M) \arrow[d] \\
    \Maps^\gen\left(X, \BB V \right) \arrow[r] & \Maps^\gen\left(X, \YY \right) \arrow[r] & \Maps^\gen(X, \BB M)
\end{tikzcd}
\]
The composite of the lower row induces injections on stabilizers of geometric points, and hence so does the composite of the middle row. Therefore, it suffices to show that the map ($\star$) has the same property. But in fact this property holds for the composite
    \[\Maps(X, \YY) \to \Maps^\gen(X, \YY)\]
because $\YY$ is the quotient of a scheme by an algebraic group.

%% file: cofiber.tex
\section{The cofiber of $\eta_G$ is Eisenstein series}
\label{section:cofib_Eis}

\subsubsection{} The purpose of this section is to prove Lemma \ref{lemma:cofib_eta_Eis}. We will use $\Whit(\Bun_N)$ and not $\Whit\left(\Bunc_N \right)$.

\subsection{A compatibility}

\subsubsection{} The purpose of this subsection is to prove Proposition \ref{proposition:poinc_identification}. Let us define the map (\ref{eqn:Poinc*!2}) first. Consider the factorization
    \[
    \Bun_N \overset{r}{\longrightarrow} \Bun_N/Z_G \overset{\pi}{\longrightarrow} \Bun_G.
    \]
The action of $Z_G$ on $\Bun_N$ is trivial, so we have
    \[r_! \, [z_G] \simeq r_*.\]
Combining this with the forgetting supports map for $\pi$ gives the desired morphism
    \[\Poinc_! \, [z_G] \to \Poinc_*.\]

\subsubsection{Notation} We write
    \[
    r_{G,*} = j_* \cdot r_* \cdot k_{\Bun_G} \simeq j_* \cdot r_* \cdot \omega_{\Bun_G} [-2 \delta_G] \in \DMod\left(\Bunc_G \right).
    \]
There is a canonical map
    \[r_{G,!} \, [z_G] \to r_{G,*}.\]

\subsubsection{Proof of Proposition \ref{proposition:poinc_identification}} Recall that (\ref{eqn:Poinc*!1}) is defined by applying the forgetting supports morphisms $\pi_{2,!} \to \pi_{2, *}$ to the composite
\begin{equation}
    \label{eqn:Poinc*!intermediate}
    \left( \WW \stoximes r_{G,!} \right)^\trans [z_G]
    \to \left( \WW \shoximes r_{G,!} \right)^\trans [2 \delta_G + z_G]
    \to \left( \WW \shoximes r_{G,*} \right)^\trans [2 \delta_G].
\end{equation}
The source of this morphism is $!$-extended from the open substack
    \[
    j : \Bun_N \times \BB Z_G \mon \left( \Bun_N \!\underset{\Bun_G}{\times} \Bunc_G \right)^\trans
    \]
and the target is $*$-extended from the same locus. Therefore, (\ref{eqn:Poinc*!1}) identifies with the forgetting supports morphism for
    \[\Bun_N \times \BB Z_G \to \Bun_G\]
applied to the restriction of (\ref{eqn:Poinc*!intermediate}) along $j$. This restriction identifies with the canonical map
    \[r_! \cdot \WW \, [z_G] \to r_* \cdot \WW.\]

\subsection{The cofiber is Eisenstein series}

\subsubsection{} Let us form the diagram below.
\[
\begin{tikzcd}[
  column sep={6em,between origins},
  row sep={4em,between origins}
]
  \Bun_N \arrow[r, "r"] \arrow[d] & \Bun_N / Z_G \arrow[r, "\rho"] \arrow[d] & \Bun_N / T \arrow[r] \arrow[d, "j"] & \Bun_B \arrow[d]  \\
  \Bunc_N \arrow[r] & \Bunc_N / Z_G \arrow[r] & \Bunc_N / T \arrow[r, "i"] & \Bunc_B
\end{tikzcd}
\]
Observe that $i$ is a closed immersion because it is a base change of $\BB T \mon \Bun_T$. So the morphism
    \[\Poinc_! [z_G] \to \Poinc_*\]
is induced by the proper pushforward of
    \[
    j_! \cdot \rho_! \cdot r_!
    \overset{\alpha}{\longrightarrow} j_* \cdot \rho_! \cdot r_!
    \overset{\beta}{\longrightarrow} j_* \cdot \rho_* \cdot r_!
    \simeq j_* \cdot \rho_* \cdot r_* \, [-z_G].
    \]
along
    \[\pi : \Bunc_N/T \to \Bun_G.\]
Therefore, Proposition \ref{proposition:poinc_principal_series} is a consequence of the following two assertions:
\begin{itemize}
    \item The cofiber of $\pi_* \cdot \beta$ belongs to the image of $\Eis_{B,*}$. This will follow from Proposition \ref{proposition:Laumonesque}.
    \item The cofiber of $\pi_* \cdot \alpha$ lies in the full subcategory generated by all $\Eis_{P,*}$. This will follow from Proposition \ref{proposition:unipotent_gerbe}.
\end{itemize}

\begin{proposition}
\label{proposition:Laumonesque}
The object
  \[
  \cofib\left( \rho_! \cdot r_! \cdot \WW \to \rho_* \cdot r_! \cdot \WW \right)
  \]
lies in the image of
  \[q^! : \DMod(\Bun_T) \to \DMod(\Bun_B).\]
\end{proposition}

\subsubsection{} 
\label{subsubsection:laumon_calculation}
The mechanism responsible for this result is the following computation of G. Laumon \cite{laumon_fourier}. Consider the diagram:
\[
\begin{tikzcd}[
  column sep={5em,between origins},
  row sep={3em,between origins}
]
    & \Ga \arrow[d, "\varrho"] & \\
  \pt \arrow[r, "j"] & \Ga/\Gm & \arrow[l, "i" above] \BB\Gm
\end{tikzcd}
\]
By the contraction principle,
    \[i^* \cdot \varrho_* \cdot \AS \simeq \Gamma(\Ga/\Gm, \varrho_* \cdot \AS) \simeq 0.\]
Then one obtains
    \[\rho_* \cdot \AS \simeq j_! \cdot k[1]\]
by using base change to compute $j^! \cdot \varrho_* \cdot \AS$. Therefore, we have fiber sequences
    \[
    k\,[1] \to \varrho_! \cdot \AS \to i_* \cdot k_{\BB\Gm}
    \qquad
    i_* \cdot k_{\BB \Gm} \to \varrho_* \cdot \AS \to k \, [1]
    \]
One computes that the forgetting supports map for $\AS$ identifies with
    \[\varrho_! \cdot \AS \to i_* \cdot k_{\BB\Gm} \to \varrho_* \cdot \AS,\]
so its cofiber is an extension of $k[1]$ by itself.

\subsubsection{Proof} Recall that the Whittaker character is defined via a composite
    \[
    \Bun_N \overset{\ev}{\longrightarrow} N^\ab \simeq T_\ad^+ \overset{\mu}{\longrightarrow} \Ga.
    \]
Since the map $\ev$ is $T$-equivariant, we may form the Cartesian diagram
\[
\begin{tikzcd}[
  column sep={8em,between origins},
  row sep={4em,between origins}
]
  \Bun_N \arrow[r, "r"] \arrow[d, "\ev"] & \Bun_N \times \BB Z_G \arrow[r, "\rho"] \arrow[d] & \Bun_N / T \arrow[d, "\ev/T"]  \\
  T_\ad^+ \arrow[r] & T_\ad^+ \times \BB Z_G \arrow[r, "\varrho \times \id"] & T_\ad^+ / T_\ad \times \BB Z_G.
\end{tikzcd}
\]
Since $\ev$ is smooth, we have
\begin{align*}
    & (\rho_! \to \rho_*) \cdot r_! \cdot \WW \\
    & \qquad \simeq (\rho_! \to \rho_*) \cdot (\ev \times \id)^! \cdot (\mu^! \cdot \AS \boxtimes r_! \cdot k) \\
    & \qquad \qquad \simeq (\ev/T)^! \cdot (\varrho_! \to \varrho_*) \cdot (\mu^! \cdot \AS \boxtimes r_! \cdot k )
\end{align*}
by base change. Therefore, it suffices to prove that the object
  \[\cofib\left(\varrho_! \cdot \mu^! \cdot \AS \to \varrho_* \cdot \mu^! \cdot \AS \right)\]
lies in the full subcategory of $\DMod(T_\ad^+/T_\ad)$ generated by the constant sheaf. But $\varrho$ identifies with
  \[\Ga \times \cdots \times \Ga \to (\Ga / \Gm) \times \cdots \times (\Ga / \Gm)\]
and
  \[\mu^! \cdot \AS \simeq \AS \boxtimes \cdots \boxtimes \AS \]
because $\AS$ is a character sheaf. Therefore, the aforementioned calculation (\S\ref{subsubsection:laumon_calculation}) induces a filtration on the cofiber in question by the power poset of $\Delta$, each of whose graded pieces is
    \[k[1] \boxtimes \cdots \boxtimes k[1].\]

\begin{proposition}
\label{proposition:unipotent_gerbe}
The direct image along the proper morphism $\Bunc_N/T \to \Bun_G$ of any object supported on the complement of
    \[\Bun_N/T \mon \Bunc_N/T\]
is Eisenstein series.
\end{proposition}

\subsubsection{} In what follows, we will prove this proposition in several steps under the assumption that $X$ does not have genus zero. In the genus zero case, the entire category $\DMod(\Bun_G)$ is Eisenstein series if $G$ is not a torus. The proposition is vacuous when $G$ is a torus.

\subsubsection{First step} To analyse this complement, recall that $\Bunc_B$ has a defect stratification indexed by dominant cocharacters
  \[
  \left(X^\lambda \times \Bun_T \right) \underset{\Bun_T}{\times} \Bun_B \simeq \Bunc_B^\lambda \mon \Bunc_B.
  \]
Here we are using the map
    \[
    X^\lambda \times \Bun_T \to \Bun_T 
    \qquad (D, \FF_T) \mapsto \FF_T(D).
    \]
We will use the notation $\Bunc_N^\lambda / T$ for the induced stratification of $\Bunc_N / T$. It suffices to prove that for every nonzero $\lambda$, the image of $*$-pushforward along the map
    \[\Bunc_N^\lambda/T \to \Bun_G\]
is Eisenstein series.

\subsubsection{Second step} Let $P$ be an arbitrary standard parabolic. Then we have the following Cartesian diagram:
\[
\begin{tikzcd}[
  column sep={5.5em,between origins},
  row sep={4em,between origins}
]
  \Bunc_B^\lambda \arrow[r] \arrow[d] & \Bun_B \arrow[r] \arrow[d] & \Bun_P \arrow[d] \\
  \Bunc_{B_M}^\lambda \arrow[r] \arrow[d] & \Bun_{B_M} \arrow[r] \arrow[d] & \Bun_M \\
  X^\lambda \times \Bun_T \arrow[r] & \Bun_T &
\end{tikzcd}
\]
This gives the Cartesian left face of the following diagram:
\[
\begin{tikzcd}[
  column sep={6em,between origins},
  row sep={2em,between origins}
]
   & \Bunc_N^\lambda / T \arrow[r] \arrow[dl, " 'q" above left] \arrow[dd] & \Bunc_N / T \arrow[ddd, bend left=10]  \\
  \Bunc_{N_M}^\lambda / T \arrow[dd] & \\
  & \Bun_P \arrow[dl] \arrow[dr] & \\
  \Bun_M & & \Bun_G
\end{tikzcd}
\]
By base change, it suffices to show that for every nonzero $\lambda$, there exists a proper parabolic for which the image of the functor
  \[
  {'q}^! : \DMod\left( \Bunc_{N_M}^\lambda / T \right) \to \DMod\left( \Bunc_N^\lambda / T \right)
  \]
generates the target.

\subsubsection{Third step} For each $\lambda$, let $P_\lambda$ be the parabolic corresponding to the set of positive roots $\alpha$ for which $\lambda(\alpha) \neq 0$. We claim that if the genus of $X$ is nonzero, then after taking $P = P_\lambda$ in the above notation, $'q$ is tower of unipotent gerbes. This will prove the proposition because $P_\lambda$ is a proper parabolic whenever $\lambda \neq 0$.

\subsubsection{Facts} Let us recall some basic facts before proving the claim.
\begin{itemize}
    \item[(i)] Let $H \to G$ be a map of groups. Then the fiber of $\Bun_H \to \Bun_G$ over a point $\FF_G$ classifies sections of \[\FF_G \overset{G}{\times} G/H.\]
    \item[(ii)] Let $\Gm$ act on $\BB \Ga$ through the dilation action on $\Ga$. For any line bundle $\LL$ on $X$, sections of $\LL \overset{\Gm}{\times} \BB \Ga$ are classified by $\Gamma(X, \LL[1])$.
\end{itemize}
Note that (ii) is obtained from (i) by taking $H \mon \GL_2$ to be the automorphism group of $\AA^1$.

\subsubsection{Notation} Choose a sequence of quotients
  \[B = B_0 \epi B_1 \epi \cdots \epi B_n = B/U_P\]
with the property that
  \[\ker(B_i \epi B_{i+1}) \simeq \Ga\]
is central in the unipotent radical of $B_i$. Then the action of $B_{i+1}$ on $B_{i+1}/B_i \simeq \BB \Ga$ factors through $B_{i+1} \epi T$. The weight of the $T$-action is some positive root $\alpha_i$ of $M$.

\subsubsection{Proof of the claim} It suffices to show that the map $\aph q_i$ in the Cartesian diagram below is a unipotent gerbe.
\[
\begin{tikzcd}[
  column sep={6em,between origins},
  row sep={4em,between origins}
]
  \YY_i \arrow[r, "\aph q_i"] \arrow[d] & \YY_{i+1} \arrow[r] \arrow[d] & X^\lambda \arrow[d, "\AJ_\omega"] \\
  \Bun_{B_i} \arrow[r] & \Bun_{B_{i+1}} \arrow[r] & \Bun_T
\end{tikzcd}
\]
So fix a point of $\YY_{i+1}$, and say that its image in $X^\lambda$ is $D$. Then the fiber of $\aph q_i$ over our chosen point is
\begin{equation}
\label{eqn:unipotent_gerbe_complex}
    \Gamma\left(X, \omega^{\alpha_i(\rhov)}(\alpha_i(D))[1] \right)
    \simeq \Gamma \left(X, \omega^{1 - \alpha_i(\rhov)} (-\alpha_i(D)) \right)^\vee
\end{equation}
We claim that this complex is concentrated in negative degrees. Indeed, the line bundle appearing in the RHS of this expression has degree
    \[
    (1 - \alpha_i(\rhov)) \cdot \deg(\omega_X) - \alpha_i(D).
    \]
This integer is negative because $\alpha_i(\rhov)$ is a positive integer, $\alpha_i(D)$ is a positive divisor of degree $\alpha_i(\lambda) > 0$, and $\deg(\omega_X) \ge 0$ by our assumption on the genus.

%% file: Jacquet.tex
\section{The Jacquet functor}
\label{section:Jacquet}

\subsubsection{} In this section we will define a certain category $\Whit(M)_\Conf$ along with functors
  \[\Jac_\Conf : \Whit(G) \to \Whit(M)_\Conf  \qquad  \Poinc_! : \Whit(M)_\Conf \to \DMod(\Bun_M).\]
Furthermore, we will construct a commutative square
\begin{equation}
\label{eqn:Jac_Conf_Poinc_compatibility}
\begin{tikzcd}[
  column sep={9em,between origins},
  row sep={4em,between origins}
]
  \Whit(G) \arrow[d, "\Poinc_!" left] \arrow[r, "\Jac_\Conf"] & \Whit(M)_\Conf \arrow[d, "\Poinc_!"] \\
  \DMod(\Bun_G) \arrow[r, "\CT_{P,!}^-"] & \DMod(\Bun_M) 
\end{tikzcd}
\end{equation}
modulo a certain ULA result whose proof we defer to the next section.

\subsection{Definition of the Jacquet functor}

\subsubsection{Zastava spaces} Write
\[
    \ZZt = \ZZ_{\Nc,\Pt^-} \mon \Bunc_N \underset{\Bun_G}{\times} \Bunt_{P^-}
\]
for the open substack where the generic $N$ and $P^-$ reductions are transverse. We will use the notation $(-)_\ZZt$ to denote restriction to $\ZZt$. We will also use notation such as $\ZZ_{\Bc, \Pt^-}$, $\ZZ_{\Nc, P^-}$, and $\ZZ_{N, P^-}$ to denote the similarly defined stacks with other subgroups and degeneracy conditions.

\subsubsection{Configuration space} Recall that we have a map
  \[\Lambda_G \to \Lambda_{G,P}\]
from the cocharacter lattice of $T$ to the cocharacter lattice of $M^\ab$. We write $\Lambda_{G,P}^\pos$ for the submonoid generated by the images of the simple roots not lying in the Dynkin subdiagram corresponding to $M$. We write
  \[\Conf_{G,P} = \sqcup_{\theta \in \Lambda_{G,P}^\pos} X^\theta\]
for the space of divisors colored by $\alpha_i$, $i \notin \Delta_M$.

\subsubsection{} Let us regard $X$ as a constant curve over $\Conf_{G,P}$ and write $\overset{\circ}{X} \mon X$ for the complement of the tautological divisor.
\begin{itemize}
    \item We define
    \begin{align*}
        \Bunc_{B_M, \Conf}^\vpol &= \Maps^\gen \left(\overset{\circ}{X}, \BB B_M \mon M \bs \left(\overline{M / N_M}\right)/T \right)  \\
        \Bunc_{B_M, \Conf}^\pol &= \Bunc_{B_M, \Conf}^\vpol \underset{\Maps(\overset{\circ}{X}, \BB M)}{\times} \left( \Bun_M \times \Conf \right).
    \end{align*}
    \item Consider the map
        \[\Conf \to \Maps\left(\overset{\circ}{X}, \BB T\right)\]
    given by restricting $\omega^\rhov$ to $\overset{\circ}{X}$. We define
    \begin{align*}
    \Bunc_{N_M, \Conf}^\vpol &= \Bunc_{B_M, \Conf}^\vpol \underset{\Maps(\overset{\circ}{X}, \BB T)}{\times} \Conf  \\
    \Bunc_{N_M, \Conf}^\pol &= \Bunc_{B_M, \Conf}^\pol \underset{\Maps(\overset{\circ}{X}, \BB T)}{\times} \Conf.
    \end{align*}
\end{itemize}
So $\Bunc_{N_M, \Conf}^\vpol$ classifies a point of $\Conf$ along with a degenerate $N_M$-bundle on $\overset{\circ}{X}$. The stack $\Bunc_{N_M, \Conf}^\vpol$ classifies this data along with an extension of the induced $M$-bundle to all of $X$. The reason for introducing $\Bunc_{N_M, \Conf}^\pol$ is that it receives a natural map from $\ZZt$. The following will appear in \S\ref{subsection:Zastava_spaces}.

\begin{construction}
\label{construction:zast_u_lift}
The map $\ZZt \to \Bun_M$ admits a lift
    \[u : \ZZt \to \Bunc_{N_M, \Conf}^\pol.\]
Furthermore, this morphism is proper.
\end{construction}

\subsubsection{} The definition of the full subcategory
    \[
    \Whit(M)_\Conf = \Whit\left(\Bunc_{N_M, \Conf}^\pol\right) \mon \DMod\left(\Bunc_{N_M, \Conf}^\pol\right).
    \]
follows a standard pattern. Write
    \[\left(\Bunc_{N_M, \Conf}^\pol \times \Ran \right)^\good \mon \Bunc_{N_M, \Conf}^\pol \times \Ran\]
for the open subset cut out by the condition that the point of $\Ran$ is disjoint from both the support of the divisor and the degeneracy locus of the $N_M$-reduction. Let
\begin{equation}
\label{eqn:Whit_Conf_Hecke}
\begin{tikzcd}[
  column sep={10em,between origins},
  row sep={2em,between origins}
]
    & \HH_{N_M, \Conf} \arrow[dl] \arrow[dr] \arrow[dd, "\psi"] & \\
    \left(\Bunc_{N_M, \Conf}^\pol \times \Ran \right)^\good & & \left(\Bunc_{N_M, \Conf}^\pol \times \Ran \right)^\good \\
    & \Ga &
\end{tikzcd}
\end{equation}
denote the Hecke groupoid which acts by modifying the $N_M$-bundle, equipped with its natural character. We will say that $\WW$ belongs to $\Whit(M)_\Conf$ if 
    \[(\WW \boxtimes \omega_{\Ran})^\good\]
is equivariant against $\psi^! \cdot \AS$.

\subsubsection{} The functor
    \[
    \oblv : \Whit(M)_\Conf \mon \DMod\left(\Bunc_{N_M, \Conf}^\pol\right)
    \]
admits a continuous right adjoint
    \[
    \Av_* : \DMod\left(\Bunc_{N_M, \Conf}^\pol\right) \to \Whit(M)_\Conf 
    \]
and we define
    \[
    \Poinc_{\Conf,!} = \pi_{\Conf,!} \cdot \oblv
    \qquad
    \coeff_{\Conf,*} = \Av_* \cdot \pi_{\Conf}^!,
    \]
where $\pi_\Conf$ is the projection to $\Bun_M$ displayed in the diagram (\ref{eqn:zastava_diagram}) below.

\subsubsection{} Now we are ready to describe the construction of $\Jac_{\Conf}$. Consider the diagram:
\begin{equation}
\label{eqn:zastava_diagram}
\begin{tikzcd}[
  column sep={4em,between origins},
  row sep={4em,between origins}
]
  & & \ZZt \arrow[dl, hook] \arrow[ddrr, dashed, "u"] & & \\
  & \Bunc_N \! \underset{\Bun_G}{\times} \! \Bunt_{P^-} \arrow[dl, "\aph \widetilde{p}"] \arrow[dr, "\pi"] & & & \\
  \Bunc_N \arrow[dr] & & \Bunt_{P^-} \arrow[dl, "\widetilde{p}"] \arrow[dr, "\qt"] & & \Bunc_{N_M, \Conf}^\pol \arrow[dl, "\pi_\Conf"] \\ 
  & \Bun_G & & \Bun_M &
\end{tikzcd}
\end{equation}
Write $j_{P,!}$ for the $!$-extension of the constant sheaf from $\Bun_{P^-}$ to $\Bunt_{P^-}$. The following assertion is standard, but we will recall its proof at the end of the section.

\begin{proposition}
\label{proposition:zast_vanishing}
Let $i$ denote the closed complement to $\ZZt \mon \Bunc_N \! \underset{\Bun_G}{\times} \! \Bunt_{P^-}$. Then the $!$-pushforward of
  \[i^* \left( \WW \stoximes j_{P,!}^- \right)\]
to $\Bun_M$ is zero.
\end{proposition}

\subsubsection{} This vanishing gives an isomorphism
    \[
    \CT_P^- \cdot \Poinc_! \cdot \WW 
    \simeq (\qt \cdot \pi)_! \cdot \left( \WW \stoximes j_{P,!}^- \right)_\ZZt.
    \]
Since $u$ is proper, the map
    \[
    (\qt \cdot \pi)_! \cdot \left( \WW \stoximes j_{P,!}^- \right)_\ZZt
    \simeq \pi_{\Conf,!} \cdot u_! \cdot \left( \WW \stoximes j_{P,!}^- \right)_\ZZt
    \to \pi_{\Conf,!} \cdot u_* \cdot \left( \WW \shoximes j_{P,!}^- \right)_\ZZt \, [2\delta_G].
    \]
is an isomorphism by the following result, whose proof will occupy \S\ref{section:ULA}.

\begin{lemma}
\label{lem:zast_ULA}
The purity morphism
    \[
    \left( \WW \stoximes j_{P,!}^- \right)_\ZZt 
    \to \left( \WW \shoximes j_{P,!}^- \right)_\ZZt \, [2 \delta_G]
    \]
is an isomorphism.
\end{lemma}

\subsubsection{} Therefore, we define
    \[\Jact_\Conf = u_* \cdot \left( - \stoximes j_{P,!}^- \right)_\ZZ \, [2 \delta_G].\]
and consider the diagram
\[
\begin{tikzcd}
  \Whit(G) \arrow[d] \arrow[r, dashed, "\Jac"] & \Whit(M)_\Conf \arrow[d] \\
  \DMod\left( \Bunc_N \right) \arrow[r, "\Jact"] & \DMod\left( \Bunc_{N_M,\Conf}^\pol \right)
\end{tikzcd}
\]
We define $\Jac_\Conf$ to be the upper dashed arrow, which makes by virtue of the proposition below. It follows that
  \[\CT_{P,!}^- \cdot \Poinc_! \simeq \Poinc_! \cdot \Jac_\Conf.\]

\begin{proposition}
\label{proposition:Jac_Whit_preservation}
The functor $\Jact_\Conf$ takes Whittaker sheaves to Whittaker sheaves.
\end{proposition}

\subsection{Some invariant theory}

\subsubsection{} Construction \ref{construction:zast_u_lift} is based on some elementary facts about the invariant theory of $G$ and its subgroups. Everything we will need is summarized in Proposition \ref{proposition:zast_cartesian}, whose proof can be obtained by a close reading of \cite[\S3]{wang_monoid}.

\subsubsection{Notation} For a classical affine scheme $X = \Spec(A)$ with an action of a group $G$, we write
    \[X \git G = \Spec(A^G)\]
for the naive GIT quotient. Here we are taking invariants in the underived sense. So, for example, we have $\overline{G/N} = G \git N$ in this notation.

\subsubsection{} Recall that the closure of
  \[M \simeq P/U_P \mon G \git U_P\]
is a reductive monoid containing $M$ as its group of units. We will denote it by $M^+$. Similarly, the closure of
  \[M^\ab \simeq P/[P,P] \mon G \git [P,P]\]
is a reductive monoid containing $M^\ab$ as its group of units. We denote it $M^{\ab,+}$. 

\begin{proposition}
\label{proposition:zast_cartesian}
Consider the commutative diagram
\begin{equation}
  \begin{tikzcd}[
  column sep={9em,between origins},
  row sep={4em,between origins}
]
\label{eqn:zast_cartesian}
    U_P^- \bss P^- \overset{M}{\times} P \git N \arrow[r, hook] \arrow[d] & U_P^- \bss G \overset{G}{\times} G \git N \arrow[d] \\
    U_P^- \bss P^- \overset{M}{\times} P \git [P,P] \arrow[r, hook] & U_P^- \bss G \git [P,P]
  \end{tikzcd}
\end{equation}
of affine schemes with $M \times T$ action.
\begin{itemize}
    \item[(i)] This diagram is Cartesian.
    \item[(ii)] The lower row identifies with $M^\ab \mon M^{\ab,+}$.
\end{itemize}
\end{proposition}

\subsubsection{Order of operations} Our convention is that GIT quotient takes precedence over balanced product.

\subsubsection{Proof of Proposition \ref{proposition:zast_cartesian}.(ii)} Actually (ii) is a consequence of a more precise assertion, which we state separately: 

\begin{proposition}
\label{proposition:retract}
The composite
  \[M^{\ab,+} \to G \git [P,P] \to U_P^- \bss G \git [P,P]\]
is an isomorphism.
\end{proposition}

\subsubsection{Proof} Consider the diagram below obtained by functoriality.
\[
\begin{tikzcd}[
  column sep={7em,between origins},
  row sep={3em,between origins}
]
  M^+ \arrow[r, hook] \arrow[d] & G \git U_P \arrow[r] \arrow[d] & U_P^- \bss G \git U_P \arrow[d] \\
  M^+ \git [M,M] \arrow[r, hook] & G \git [P,P] \arrow[r] & U_P \bss G \git [P,P]
\end{tikzcd}
\]
According to \cite[Theorem 3.2.8]{wang_monoid}, the upper row is a retract. Therefore, the lower row is also a retract. In fact, it is a retract of affine schemes, so the lower left arrow is a closed immersion. Therefore, $M^+ \git [M,M]$ contains $M^{\ab,+}$ as a closed subscheme. But $M$ maps dominantly onto both, so they are equal.

\subsubsection{} In order to prove Proposition \ref{proposition:zast_cartesian}.(i), it will be useful to have some generalities on GIT quotients. Let $\UU \mon X$ be integral affine schemes with $G$-action. We will say that the pair $(\UU \mon X, G)$ is \emph{good} if the following square is Cartesian:
  \[
  \begin{tikzcd}[
  column sep={4em,between origins},
  row sep={3em,between origins}
]
    \UU \arrow[r] \arrow[d] & X \arrow[d] \\
    \UU\git G \arrow[r] & X \git G
  \end{tikzcd}
  \]
We will further say that the pair is \emph{quite good} if the lower arrow is a distinguished open affine.

\subsubsection{Example} Let $\UU \mon \AA^2$ be the complement of a line through the origin, and let $\Gm$ act by dilation. Then the pair $(\UU \mon X, \Gm)$ is not good.

\begin{proposition}
\label{proposition:git_qgood}
Let $X$ be an integral affine scheme with $G$-action.
\begin{itemize}
    \item[(i)] For any distinguished open affine $V \mon X \git G$, the pair
      \[\left( V \! \underset{X \git G}{\times} \! X \mon X, G \right)\]
    is quite good.
    \item[(ii)] Suppose that $(\UU \mon X, G)$ is quite good. For any subgroup $H \mon G$, the pair $(\UU \mon X, H)$ is quite good and the square
      \[
      \begin{tikzcd}[
  column sep={4em,between origins},
  row sep={3em,between origins}
]
        \UU \git H \arrow[r, hook] \arrow[d] & X \git H \arrow[d] \\
        \UU \git G \arrow[r, hook] & X \git G
      \end{tikzcd}
      \]
    is Cartesian.
    \item[(iii)] Let $H \mon G$ be a normal subgroup. If $(\UU \mon X, H)$ is good and $(\UU \git H \mon X \git H, G/H)$ is quite good, then $(\UU \mon X, G)$ is quite good.
\end{itemize}
\end{proposition}

\subsubsection{Proof of (i)} Let $\UU$ denote the fiber product indicated by the diagram.
\[
\begin{tikzcd}[
  column sep={3.5em,between origins},
  row sep={3em,between origins}
]
   & \UU \arrow[dl] \arrow[d] \arrow[r, hook] & X \arrow[d] \\
  \UU \git G \arrow[r, dashed] & V \arrow[r, hook] & X \git G
\end{tikzcd}
\]
We need to check that the dashed arrow is an isomorphism. Write $X = \Spec(A)$ and let $f \in A^G$ be a function cutting out the complement of $V$. Then the dashed arrow corresponds to the inclusion of rings
    \[A^G[\tfrac{1}{f}] \mon A[\tfrac{1}{f}]^G.\]
This is an isomorphism.

\subsubsection{Proof of (ii)} Write $V$ for the fiber product appearing in the following Cartesian diagram:
\[
\begin{tikzcd}[
  column sep={4em,between origins},
  row sep={3em,between origins}
]
   & \UU \arrow[dl] \arrow[d] \arrow[r, hook] & X \arrow[d] \\
  \UU \git H \arrow[r, dashed] & V \arrow[r, hook] \arrow[d] & X \git H \arrow[d] \\
   & \UU \git G \arrow[r, hook] & X \git G
\end{tikzcd}
\]
Note that $V \mon X \git H$ is a distinguished open affine because $\UU \git G \mon X \git G$ is. Therefore we may apply (i).

\subsubsection{Proof of (iii)} This is a formal consequence of the pullback lemma.

\subsubsection{Proof of Proposition \ref{proposition:zast_cartesian}.(i)} We give a sequence of reductions.

\subsubsection{First reduction} It suffices to show that (\ref{eqn:zast_cartesian}) is Cartesian after pullback along the $G$-torsor
  \[U_P^- \bss G \times G \git N \to U_P^- \bss G \overset{G}{\times} G \git N.\]
So we need to show that the diagram
\[
\begin{tikzcd}[
  column sep={10em,between origins},
  row sep={4em,between origins}
]
  U_P^- \bss P^- \overset{M}{\times} P \times G \git N \arrow[r, hook] \arrow[d] & U_P^- \bss G \times G \git N \arrow[d] \\
  U_P^- \bss P^- \overset{M}{\times} P \git [P,P] \arrow[r, hook] & U_P^- \bss G \git [P,P]
\end{tikzcd}
\]
is Cartesian. Let us spell out the precise definition of the upper arrow. Consider the pair
\begin{equation}
\label{eqn:git_pair}
    \left(P^- \overset{M}{\times} P \times G \mon G \times G, U_P^- \times G \times [P,P] \right)    
\end{equation}
where the map is
    \[(x, h) \mapsto (x h^{-1}, h).\]
and the group actions are
    \[
    (u, g, p) \cdot (x, h) = (u \cdot x \cdot p, g \cdot h \cdot p)
    \qquad
    (u, g, p) \cdot (h, h') = (u \cdot h \cdot g^{-1}, g \cdot h' \cdot p).
    \]
Taking the GIT quotient by the subgroup $U_P^- \times N \mon U_P^- \times G \times [P,P]$ gives the upper arrow.

\subsubsection{Second reduction} Now Proposition \ref{proposition:git_qgood}.(ii) implies that it suffices to show that (\ref{eqn:git_pair}) is quite good. Note that
    \[
    \left(P^- \overset{M}{\times} P \times G \mon G \times G, G \right)
    \]
is good. Therefore, by Proposition \ref{proposition:git_qgood}(iii), it is enough to prove that
    \[
    \left(P^- \overset{M}{\times} P \mon G, U_P^- \times [P,P] \right)
    \]
is quite good. 

\subsubsection{Third reduction} By Proposition \ref{proposition:retract}, the diagram
\[
\begin{tikzcd}[
  column sep={12em,between origins},
  row sep={5em,between origins}
]
  U_P^- \bss \left( P^- \overset{M}{\times} P \right) \git U_P \arrow[r, hook] \arrow[d] & U_P^- \bss G \git P \arrow[d] \\
  U_P^- \bss \left( P^- \overset{M}{\times} P \right) \git [P,P] \arrow[r, hook] & U_P^- \bss G \git [P,P]
\end{tikzcd}
\]
identifies with
\[
\begin{tikzcd}[
  column sep={5em,between origins},
  row sep={4em,between origins}
]
  M \arrow[r, hook] \arrow[d] & M^+ \arrow[d] \\
  M^\ab \arrow[r, hook] & M^{\ab, +}
\end{tikzcd}
\]
This square is Cartesian by \cite[Appendix C.4.3]{lin_nearby_cycles}. Therefore, by Proposition \ref{proposition:git_qgood}.(iii), it suffices to show that
  \[\left(P^- \overset{M}{\times} P \mon G, U_P^- \times U_P\right)\]
is good. 

\subsubsection{End of the proof} It remains to show that the diagram
\begin{equation}
\label{eqn:git_cartesian}
\begin{tikzcd}[
  column sep={5em,between origins},
  row sep={4em,between origins}
]
  P^- \overset{M}{\times} P \arrow[r, hook] \arrow[d] & G \arrow[d] \\
  M \arrow[r, hook] & M^+ 
\end{tikzcd}
\end{equation}
is Cartesian. For $P = B$ one can verify this by direct calculation \cite[\emph{Proof} of Lemma 3.2.10]{wang_monoid}. The corresponding diagram is displayed on the left.
\[
\begin{tikzcd}[
  column sep={5em,between origins},
  row sep={4em,between origins}
]
  B^- \overset{T}{\times} B \arrow[r, hook] \arrow[d] & G \arrow[d] \\
  T \arrow[r, hook] & N^- \bss G \git N 
\end{tikzcd}
\qquad
\begin{tikzcd}[
  column sep={5em,between origins},
  row sep={4em,between origins}
]
  B^- \overset{T}{\times} B \arrow[r, hook] \arrow[d] & G \arrow[d] \\
  B_M^- \overset{T}{\times} B_M \arrow[r, hook] & M^+ 
\end{tikzcd}
\]
Proposition \ref{proposition:git_qgood}.(ii) now implies that the diagram on the right is also Cartesian. But this diagram will remain Cartesian if we replace the left vertical arrow with its $M$-orbit. The result of this operation is (\ref{eqn:git_cartesian}).

\subsection{Zastava spaces}
\label{subsection:Zastava_spaces}

\subsubsection{} Now we will obtain a precise description of $\ZZt$ in local terms. In order to state the result, we will need to introduce some more stacks over $\Conf_{G,P}$. 

\subsubsection{} As before, let us view $X$ as a constant curve over $\Conf$. Write
    \[\Disco \mon \Disc \mon X\]
for the punctured formal disc inside the formal disc around the support of the tautological divisor. For a map $\YY_1 \to \YY_2$, we use the notation
    \[
    \Maps^\gen(\Disc, \YY_1 \to \YY_2) = \Maps(\Disc, \YY_2) \underset{\Maps(\Disco, \YY_2)}{\times} \Maps(\Disco, \YY_1).
    \]

\subsubsection{Notation}
\begin{itemize}
    \item We define
        \begin{align*}
            \Gr_{G,\Conf} &= \Maps^\gen \left( \Disc, \pt \to \BB G \right)   \\
            \Grc_{B,\Conf} &= \Maps^\gen \left( \Disc, \pt \to T \bs (N \bss G) / G \right)   \\
            \Grt_{P^-,\Conf} &= \Maps^\gen \left( \Disc, \pt \to G \bs (G \git U_P^- ) / M \right)
        \end{align*}
    \item Note that there is a tautological section
        \begin{equation}
            \label{eqn:tautological_section}
            \Conf \to \Maps^\gen(\Disc, M^\ab \bs M^{\ab, +}).
        \end{equation}
    The map
        \[N \bss G \overset{G}{\times} G \git U_P^- \to M^{\ab, +}\]
    induces a morphism
        \[\Grc_{B,\Conf} \underset{\Gr_{G,\Conf}}{\times} \Grt_{P^-, \Conf} \to \Maps^\gen(\Disc, M^\ab \bs M^{\ab, +})\]
    and we define $\ZZ_{\Bc, \Pt, \Conf}$ to be its base change along (\ref{eqn:tautological_section}).
\end{itemize}
\subsubsection{Variants}
\begin{itemize}
    \item Consider the map
        \[
        \Conf \to \Maps\left(\overset{\circ}{X}, \BB T \right)
        \]
    given by restricting $\omega^\rhov$ to $\overset{\circ}{X}$. We define
        \[
        \Bunc_{N_M, \Conf}^\vpol = \Bunc_{B_M, \Conf}^\vpol \underset{\Maps(\overset{\circ}{X}, \BB T)}{\times} \Conf.
        \]
    \item Consider the map
        \[
        \Conf \to \Maps(\Disc, \BB T)
        \]
    given by restricting $\omega^\rhov$ to $\DD$. We define
        \begin{align*}
            \Grc_{N,\Conf} &= \Grc_{B,\Conf} \underset{\Maps(\Disc, \BB T)}{\times} \{\omega_\Disc^\rhov\} \\
            \ZZ_{\Nc, \Pt^-, \Conf} &= \ZZ_{\Bc, \Pt^-, \Conf} \underset{\Maps(\Disc, \BB T)}{\times} \{\omega_\Disc^\rhov\}
        \end{align*}
\end{itemize}

\begin{proposition}
\label{proposition:zast_loc_formula}
Define $\ZZ_{\Bc, \Pt^-}$ by replacing $N$ with $B$ in the definition of $\ZZt = \ZZ_{\Nc, \Pt^-}$. There are canonical isomorphisms
    \begin{align*}
        \ZZ_{\Bc, \Pt^-} & \simeq \Bunc_{B_M, \Conf}^{\vpol} \timest \ZZ_{\Bc, \Pt^-, \Conf} \\
        \ZZ_{\Nc, \Pt^-} & \simeq \Bunc_{N_M, \Conf}^{\vpol} \timest \ZZ_{\Nc, \Pt^-, \Conf}.
    \end{align*}
\end{proposition}

\subsubsection{} Let us explain the meaning of these twisted products, which \emph{prima facie} involve fiber products over $\BB \LL_\Conf B_M$ and $\BB \LL_\Conf N_M^\omega$. To avoid the problem of taking quotients by ind-infinite dimensional group schemes, we use
    \[\Maps(\Disco, \BB B_M)\]
as a substitute for $\BB \LL_\Conf B_M$. Then the twisted product involving $B$ is defined to be the pullback of 
\begin{equation}
\label{eqn:Zast_descent}
    \Maps^\gen\left(\Disc, \BB B_M \mon T \bs \left( N \bss G \overset{G}{\times} G \git U_P^- \right) / M \right)
\end{equation}
along the natural map
    \[
    \alpha : \Bunc_{B_M, \Conf}^\vpol \to \Maps(\Disco, \BB B_M).
    \]
The second twisted product involving $N$ is defined similarly.

\subsubsection{Sanity check} The expression in (\ref{eqn:Zast_descent}) is a substitute for
    \[\LL_\Conf B_M \bs \ZZ_{\Bc, \Pt^-, \Conf}\]
in the sense that its base change along
    \[\beta : \Conf \to \Maps(\Disco, \BB B_M)\]
identifies with $\ZZ_{\Bc, \Pt^-, \Conf}$. Since $\alpha$ Zariski locally factors through $\beta$, the twisted product is Zariski locally on $\Bunc_{B_M,\Conf}^\vpol$ a product
    \[\Bunc_{B_M,\Conf}^\vpol \underset{\Conf}{\times} \ZZ_{\Bc, \Pt^-, \Conf}\]
over the configuration space.

\subsubsection{Completion of Construction \ref{construction:zast_u_lift}} Using the identification of Proposition \ref{proposition:zast_loc_formula}, we define the morphism $u$ to be
    \[
    \Bunc_{N_M, \Conf}^{\vpol} \timest \ZZ_{\Nc, \Pt^-, \Conf}
    \to \Bunc_{N_M, \Conf}^{\vpol} \timest \Gr_{M,\Conf}
    \simeq \Bunc_{N_M, \Conf}^\pol.
    \]
It is proper because $\ZZ_{\Nc, \Pt^-, \Conf} \to \Conf$ is proper.

\subsubsection{Proof of Proposition \ref{proposition:zast_loc_formula}} Let us rewrite (\ref{eqn:zast_cartesian}) as
\[
  \begin{tikzcd}[
  column sep={7em,between origins},
  row sep={4em,between origins}
]
    M \git N_M \arrow[r, hook] \arrow[d] & U_P^- \bss G \overset{G}{\times} G \git N \arrow[d] \\
    M^\ab \arrow[r, hook] & M^{\ab,+}
  \end{tikzcd}
\]
The quotient of the lower row by $M \times T$ identifies with
    \[\BB [M,M] \times \BB T \mon \BB [M,M] \times M^\ab \bs M^{\ab,+} \times \BB T.\]
Therefore, we obtain a Cartesian diagram
\[
  \begin{tikzcd}[
  column sep={12em,between origins},
  row sep={4em,between origins}
]
    M \bs \left( M \git N_M \right) / T \arrow[r, hook] \arrow[d] & M \bs \left( U_P^- \bss G \overset{G}{\times} G \git N \right) / T \arrow[d] \\
    \pt \arrow[r, hook] & M^\ab \bs M^{\ab,+}
  \end{tikzcd}
\]
Now Proposition \ref{proposition:zast_loc_formula} is a formal consequence of Beauville--Laszlo gluing.

\subsection{The Whittaker condition}
\label{subsection:Jacquet_Whittaker}

\subsubsection{Proof of Proposition \ref{proposition:zast_vanishing}} We follow the usual pattern of passing to strata and applying character orthogonality.

\subsubsection{First reduction} The object $\WW \stoximes j_{P,!}$ in question is $!$ extended from the locus
  \[\Bun_N \underset{\Bun_G}{\times} \Bun_{P^-}.\]
Therefore, it suffices to show that for each locally closed substack
  \[i_w : \ZZ_{N,P^-}^w = \left( \Bun_N \underset{\Bun_G}{\times} \Bun_{P^-} \right)^w \mon \Bun_N \underset{\Bun_G}{\times} \Bun_{P^-}\]
with $w \neq 1$, the object
  \[i_w^* \left( \WW \stoximes k_{\Bun_{P^-}} \right) \simeq \aph p^* \cdot \WW \]
is annihilated by $!$-pushforward to $\Bun_M$. 

\subsubsection{Second reduction} For any finite set $x \mon X$ we may consider the open substack
  \[\ZZ_{N,P^-}^{w, \good\,x} \mon \ZZ_{N,P^-}^w\]
where the $N$ and $P^-$ reductions are in relative position $w$ at $x$. Since $\ZZ_{N,P^-}^w$ is a colimit of open substacks of this form, it suffices to show that the $*$-pullback of $\WW$ to each $\ZZ_{N,P^-}^{w, \good\,x}$ is annihilated by $!$-pushforward to $\Bun_M$. 

\subsubsection{End of the proof} Now consider the groupoid acting on $\ZZ_{N,P^-}^{w, \good \, x}$ by $N^w \cap P^-$-Hecke modifications at $x$. The unipotent subgroupoid corresponding to
    \[N^w \cap U_P^- \simeq \ker\left( N^w \cap P^- \to M \right)\]
acts along the fibers of the projection to $\Bun_M$. Now the $*$-pullback of $\WW$ is equivariant for the character associated to the pullback of $\psi$ along
  \[N^w \cap U_P^- \mon N^w \simeq N \to \Ga,\]
which is nontrivial because $w \neq 1$.

\subsubsection{Proof of Proposition \ref{proposition:Jac_Whit_preservation}} This claim is verified by doing some bookkeeping with Hecke groupoids.

\subsubsection{Step one} We form an extension of the diagram (\ref{eqn:Whit_Conf_Hecke}) appearing in the definition of Whittaker condition for $\Whit(M)_\Conf$.
\[
\begin{tikzcd}[
  column sep={10em,between origins},
  row sep={2em,between origins}
]
    & \HH_\ZZt \arrow[dl] \arrow[dr] \arrow[dd] & \\
    \left(\ZZt \times \Ran \right)^\good \arrow[dd] & & \left(\ZZt \times \Ran \right)^\good \arrow[dd] \\
    & \HH_{N_M, \Conf} \arrow[dl] \arrow[dr] \arrow[dd, "\psi"] & \\
    \left(\Bunc_{N_M, \Conf}^\pol \times \Ran \right)^\good & & \left(\Bunc_{N_M, \Conf}^\pol \times \Ran \right)^\good \\
    & \Ga &
\end{tikzcd}
\]
Here are the definitions:
\begin{itemize}
    \item We write $(\ZZt \times \Ran)^\good$ for the preimage of
        \[
        \left(\Bunc_{N_M, \Conf}^\pol \times \Ran \right)^\good \mon \Bunc_{N_M, \Conf}^\pol \times \Ran
        \]
    along the projection
        \[u \times \id : \ZZt \times \Ran \to \Bunc_{N_M,\Conf}^\pol \times \Ran.\]
    \item The fiber products of the two faces are naturally identified. We write $\HH_\ZZt$ for this common value; it is naturally a groupoid acting on $(\ZZt \times \Ran)^\good$.
\end{itemize}
Therefore, it remains to show that
    \[\left(\WW \shoximes j_{P,!}^-\right)_\ZZt\]
is $(\HH_\ZZt, \psi)$-equivariant, since $u_*$ sends such objects into $\Whit(M)_\Conf$.

\subsubsection{Step two} The $!$-pullback of $\WW$ along $\ZZt \to \Bunc_N$ is $(\HH_\ZZt, \psi)$-equivariant. Indeed, we have a commutative diagram
\[
\begin{tikzcd}[
  column sep={8em,between origins},
  row sep={2em,between origins}
]
    & \HH_\ZZt \arrow[dl] \arrow[dr] \arrow[dd, dashed] & \\
    \left(\ZZt \times \Ran \right)^\good \arrow[dd] & & \left(\ZZt \times \Ran \right)^\good \arrow[dd] \\
    & \HH_N \arrow[dl] \arrow[dr] \arrow[dd, "\psi"] & \\
    \left(\Bunc_N \times \Ran \right)^\good & & \left(\Bunc_N \times \Ran \right)^\good \\
    & \Ga &
\end{tikzcd}
\]
in which the two faces are not Cartesian, and the pullback of the Whittaker character from $\HH_N$ agrees with the pullback of the Whittaker character from $\HH_{N_M,\Conf}$.

\subsubsection{Step three} The $!$-pullback of $j_{P,!}$ along $\ZZt \to \Bunt_{P^-}$ is $\HH_\ZZt$-equivariant. Indeed, we have a commutative diagram
\[
\begin{tikzcd}[
  column sep={8em,between origins},
  row sep={2em,between origins}
]
    & \HH_\ZZt \arrow[dl] \arrow[dr] \arrow[dd, dashed] & \\
    \left(\ZZt \times \Ran \right)^\good \arrow[dd] & & \left(\ZZt \times \Ran \right)^\good \arrow[dd] \\
    & \HH_{P^-} \arrow[dl] \arrow[dr] & \\
    \left(\Bunt_{P^-} \times \Ran \right)^\good & & \left(\Bunt_{P^-} \times \Ran \right)^\good
\end{tikzcd}
\]
in which the two faces are not Cartesian, and $j_{P,!}^-$ is equivariant for $\HH_{P^-}$.

%% file: ULA.tex
\section{ULA proofs}
\label{section:ULA}

\subsubsection{} The purpose of this section is to prove Lemma \ref{lem:zast_ULA}. In the case $P = B$, this follows from the work of J. Campbell \cite{justin_Kontsevich}. The main input to that work was the Kontsevich compactification, whose parabolic version $\Bunc_P^K \to \Bunc_P$ does not seem to lift along the map $\Bunt_P \to \Bunc_P$.

\subsection{Outline of the proof}

\subsubsection{} Let us write $\ZZ \mon \ZZt$ for $\ZZ_{N,P^-} \mon \ZZ_{\Nc, \Pt^-}$. We first reformulate the assertion of the lemma in terms that are smooth-local with respect to the data of
  \[j : \ZZ  \mon  \ZZt  \qquad  \qquad \FF = \left( \WW \shoximes j_{P,!}^- \right)_\ZZt.\]
So let us examine the purity morphism
    \[
    \left( \WW \stoximes j_{P,!}^- \right)_\ZZt 
    \to \left( \WW \shoximes j_{P,!}^- \right)_\ZZt \, [2 \delta_G]
    \]
in question. Note that the LHS is $!$-extended along $j$. Furthermore, the restriction of the purity map to $\ZZ$ is an isomorphism because $\WW \boxtimes j_{P,!}$ is a local system on $\Bun_N \times \Bun_{P^-}$. Therefore, our purity morphism identifies with the canonical morphism
\begin{equation}
\label{eqn:zast_ULA_rewrite}
    j_! \cdot j^! \cdot \FF \to \FF
\end{equation}
shifted by $[2 \delta_G]$. So Lemma \ref{lem:zast_ULA} is equivalent to the following assertion, which is smooth-local in the desired sense:

\begin{lemma}
\label{lemma:zast_!_extended}
The sheaf $\FF$ is $!$-extended from $\ZZ \mon \ZZt$.
\end{lemma}

\subsubsection{Proof of Lemma \ref{lemma:zast_!_extended}, modulo the rest of the section} Write $\UU \mon \ZZt$ for the largest open subset over which (\ref{eqn:zast_ULA_rewrite}) is an isomorphism. In Construction \ref{construction:zast_correspondence}, we will form a correspondence
\[
\begin{tikzcd}[
  column sep={4em,between origins},
  row sep={2em,between origins}
]
    & \HHt \arrow[dl, "\mu" above left] \arrow[dr, "\pi"] & \\
    \ZZt & & \ZZt
\end{tikzcd}
\]
whose action preserves $\UU$ in the sense that
    \[\mu^{-1} (\UU) = \pi^{-1} (\UU).\]
This last assertion is the content of Proposition \ref{proposition:zast_corr}. 

\subsubsection{Informal description of $\HHt$} This paragraph has no mathematical meaning. According to Proposition \ref{proposition:zast_loc_formula}, a point of $\ZZt$ over $D \in \Conf_{G,P}$ may be viewed as a degenerate $N_M$-bundle away from $D$ glued to some Zastava data on the formal neighborhood of $|D|$. Actually there is a constraint on the nature of this local Zastava data in terms of the value of $D$, but we will ignore it. In these terms, the action of $\HHt$ from right to left admits the following description:
\begin{itemize}
    \item Choose a divisor $D_\good$ whose support is disjoint from that of $D$, and impose the open condition that the $N_M$-bundle is nondegenerate at $D_\good$.
    \item Cut out the $N_M$-bundle at $D_\good$ and glue in some Zastava data coming from
        \[\Gr_N \cap \Gr_{P^-} \to \Gr_N \to \Gr_{N_M}.\]
    The result is a new point of $\ZZt$ lying over the divisor $D + D_\good$.
\end{itemize}

\subsubsection{} We will then prove Proposition \ref{proposition:SS_bound}, which gives an upper bound on the singular support of $j_{P,!}^-$. This will allow us to exhibit an open substack
    \[\ZZt_\Sing \mon \UU\]
with a fairly simple membership criterion. Finally, we prove Proposition \ref{proposition:Hecke_transitivity}, which says that that every point of $\ZZt$ can be related by the correspondence $\HHt$ to a point in $\ZZt_\Sing$.

\subsection{Construction of the correspondence}

\begin{construction} 
\label{construction:zast_correspondence}
There is a correspondence
\[
\begin{tikzcd}[
  column sep={4em,between origins},
  row sep={2em,between origins}
]
   & \HH \arrow[dd, hook] \arrow[dl] \arrow[dr] & \\
  \ZZ \arrow[dd, hook] &  & \ZZ \arrow[dd, hook] \\
   & \HHt \arrow[dl, "\mu" above] \arrow[dr, "\pi" above] & \\
  \ZZt & & \ZZt
\end{tikzcd}  
\]
with the feature that both squares are Cartesian.
\end{construction}

\begin{proposition}
\label{proposition:zast_corr}
The following hold:
\begin{itemize}
    \item[(i)] The map $\mu$ is etale. 
    \item[(ii)] The map $\pi$ is smooth-locally a product.
    \item[(ii)] The two pullbacks
      \[\mu^! \cdot \left(\WW \shoximes j_{P,!}^- \right)_\ZZt \qquad \pi^! \cdot \left(\WW \shoximes j_{P,!}^- \right)_\ZZt\]
    are smooth-locally isomorphic up to tensoring with the shift of a rank one local system.
\end{itemize}
\end{proposition}

\subsubsection{Notation} Our construction of $\HHt$ will be given in terms of the local description of $\ZZt$ recorded as Proposition \ref{proposition:zast_loc_formula}. We write
    \[
    \ZZ_\Conf = \ZZ_{N, P^-, \Conf}
    \qquad \ZZt_\Conf = \ZZ_{\Nc, \Pt^-, \Conf}.
    \]
The definition of $\ZZ_{N, P^-, \Conf}$ repeats the definition of $\ZZ_{\Nc, \Pt^-, \Conf}$ given in \S\ref{subsection:Zastava_spaces} but with $\Gr_N$ and $\Gr_{P^-}$ replacing $\Grc_N$ and $\Grt_{P^-}$. There is a natural map
    \[\ZZ_\Conf \to \Gr_{N,\Conf} \to \Gr_{N_M, \Conf}.\]

\subsubsection{} We take
  \[
  \HHt = \left\{ \Bunc_{N_M, \Conf \times \Conf}^{\vpol} \right\}_\disj \underset{[\Conf \times \Conf]_\disj}{\timest} \left\{ \ZZt_\Conf \times \ZZ_\Conf \right\}_\disj
  \]
Then the leftward map is
\begin{align*}
    \mu : \HHt \mon
    & \left\{ \Bunc_{N_M, \Conf \times \Conf}^{\vpol} \right\}_\disj \underset{[\Conf \times \Conf]_\disj}{\timest} \left\{ \ZZt_\Conf \times \ZZt_\Conf \right\}_\disj \\
    & \quad \simeq \ZZt \underset{\Conf}{\times} [\Conf \times \Conf]_\disj \\
    & \qquad \to \ZZt,
\end{align*}
where we are using the sum map
    \[
    m : [\Conf \times \Conf]_\disj \to \Conf
    \]
The rightward map is
\begin{align*}
   \pi : \HHt \to 
    & \left\{ \Bunc_{N_M, \Conf \times \Conf}^{\vpol} \right\}_\disj \underset{[\Conf \times \Conf]_\disj}{\timest} \left\{ \ZZt_\Conf \times \Gr_{N_M, \Conf} \right\}_\disj \\
    & \quad \simeq \left( \ZZt \underset{\Conf}{\times} [\Conf \times \Conf]_\disj \right)^{\text{good at $D_2$}} \\
    & \qquad \mon \ZZt \underset{\Conf}{\times} [\Conf \times \Conf]_\disj \\
    & \qquad \qquad \to \ZZt
\end{align*}
where we use the map
  \[\pr_1 : [\Conf \times \Conf]_\disj \to \Conf.\]
The base change of $\ZZ \mon \ZZt$ along either map identifies with
  \[\HH = \left\{ \Bunc_{N_M, \Conf \times \Conf}^{\vpol} \right\}_\disj \underset{[\Conf \times \Conf]_\disj}{\timest} \left\{ \ZZ_\Conf \times \ZZ_\Conf \right\}_\disj.\]

\subsubsection{} This completes Construction \ref{construction:zast_correspondence}. We now turn to the proof of Proposition \ref{proposition:zast_corr}.

\subsubsection{Proofs of (i) \& (ii)} These assertions are easily verified by examining the definitions. Indeed, $\mu$ is an open immersion followed by a base change of $m$, which is etale. This proves (i). To prove (ii), we note that twisted products are smooth-locally products and that $\pr_1$ is an open subset of a product.

\subsubsection{Proof of (iii)} This argument is very similar to the proof of Proposition \ref{proposition:Jac_Whit_preservation} appearing in \S\ref{subsection:Jacquet_Whittaker}. As before, we define $(\ZZt \times \Conf)^\good$ so that the lower square of the diagram
\[
\begin{tikzcd}[
  column sep={7.5em,between origins},
  row sep={4em,between origins}
]
  & \HHt \arrow[dl, dashed, "\text{$\widehat{\mu}$ or $\widehat{\pi}$}" above left] \arrow[d, "\text{$\mu$ or $\pi$}"] \\ 
  \left(\ZZt \times \Conf\right)^\good \arrow[r] \arrow[d] & \ZZt \arrow[d] \\
  \left( \Bunc_N \times \Conf \right)^\good \arrow[r] & \Bunc_N
\end{tikzcd}
\]
is Cartesian. As the dashed arrow indicates, both $\mu$ and $\pi$ have natural lifts fitting into a commutative diagram
\[
\begin{tikzcd}[
  column sep={8em,between origins},
  row sep={2em,between origins}
]
    & \HHt \arrow[dl, "\widehat{\mu}" above left] \arrow[dr, "\widehat{\pi}"] \arrow[dd, dashed] & \\
    \left(\ZZt \times \Conf \right)^\good \arrow[dd] & & \left(\ZZt \times \Conf \right)^\good \arrow[dd] \\
    & \HH_N \arrow[dl] \arrow[dr] \arrow[dd, "\psi"] & \\
    \left(\Bunc_N \times \Conf \right)^\good & & \left(\Bunc_N \times \Conf \right)^\good \\
    & \Ga &
\end{tikzcd}
\]
Therefore, the two $!$-pullbacks of $\WW$ to $\ZZt$ differ by the $!$-pullback of the Whittaker character along $\HHt \to \HH_N$. A similar argument shows that the two $!$-pullbacks of $j_{P,!}^-$ to $\ZZt$ are isomorphic.

\subsection{A singular support estimate}

\subsubsection{} The goal of this subsection is to obtain an upper bound on the singular support of $j_{P,!}^-$. In order to state the result, we will need to describe a stratification of $\Bunt_{P^-}$ into smooth stacks
  \[{}_{\vartheta' \le \vartheta}\Bunt_{P^-} \mon \Bunt_{P^-}\]
indexed by pairs of partitions.

\subsubsection{} Recall that $\Bunc_{P^-}$ has a smooth stratification by
  \[X^\theta \times \Bun_{P^-} \simeq {}_\theta\Bunc_{P^-} \mon \Bunc_{P^-}.\]
Under this identification. the map ${}_\theta\Bunc_{P^-} \to \Bun_G$ is given by projection on $\Bun_{P^-}$, followed by the usual map $\Bun_{P^-} \to \Bun_G$.

\subsubsection{} Taking the pullback along $\Bunt_{P^-} \to \Bunc_{P^-}$, we obtain a stratification by
    \[\Bun_{P^-} \underset{\Bun_M}{\times} \HH_M^{+, \theta} \simeq  {}_\theta\Bunt_{P^-} \mon \Bunt_{P^-}.\]
These strata are not smooth because $\HH_M^{+, \theta}$ is not smooth over $\Bun_M$, so let us examine $\HH_M^{+, \theta}$ more closely. Consider first the map
    \[\HH_M^{+, \theta} \to X^\theta.\]
We have a stratification $X^\vartheta \mon X^\theta$ indexed by partitions $\vartheta \vdash \theta$. The restriction $\HH_M^{+,\vartheta}$ of $\HH_M^{+, \theta}$ to each $X^\vartheta$ is a closed substack of $\HH_{M, X^\vartheta}$. Therefore, $\HH_M^{+, \vartheta}$ inherits a relative position stratification, indexed by $\vartheta' \le \vartheta$. These strata $\HH_M^{+, \vartheta' \le \vartheta}$ are smooth over $\Bun_M \times X^\vartheta$, and we define
    \[
    {}_{\vartheta' \le \vartheta} \Bunt_{P^-} = \Bun_{P^-} \underset{\Bun_M}{\times} \HH_M^{+, \vartheta \le \vartheta}.
    \]

\begin{proposition}
\label{proposition:SS_bound}
The singular support of $j_{P,!}^-$ is contained in the conormals of the strata
    \[
    {}_{\vartheta' \le \vartheta} \Bunt_{P^-}
    \mon \Bunt_{P^-}.
    \]
\end{proposition}

\subsubsection{} The proof of Proposition \ref{proposition:SS_bound} will use Zastava spaces as a local model for $\Bunt_{P^-}$, following \cite{BFGM}. Let us recall this. For the sake of notational consistency, we write
    \[
    \Mod_M^+ = \pt \underset{\Bun_M}{\times} \HH_M^+.
    \]
There are natural maps
    \[
    \ZZ_{U_P, \Pt^-} \to \Mod_M^+ \to \Conf_{G,P}.
    \]
and the first arrow has a canonical section. We write $(-)^\theta$ to denote base change to the component $\Conf^\theta \mon \Conf$. Now according to \cite[(1)]{BFGM_erratum}, the diagram of Zastava spaces on the left
\[
\begin{tikzcd}[
  column sep={1.5em},
  row sep={1em}
]
  \Mod_M^{+,\theta} \arrow[r, hook] \arrow[d, equal] & \ZZ_{U_P, \Pt^-}^\theta & \ZZ_{U_P, P^-}^\theta \arrow[l, hook', "j_\ZZ"] \\
  \Mod_M^{+,\theta} & & 
\end{tikzcd}
\qquad 
\begin{tikzcd}[
  column sep={1.5em},
  row sep={1em}
]
  {}_{\theta}\Bunt_{P^-} \arrow[r, hook] \arrow[d] & \Bunt_{P^-} & \Bun_{P^-} \arrow[l, hook'] \\
  \HH_M^{+,\theta} & & 
\end{tikzcd}
\]
is a smooth-local model for the diagram of stacks on the right. 
\subsubsection{Proof of Proposition \ref{proposition:SS_bound}} We make a preliminary construction and the perform a sequence of reductions.

\subsubsection{A group action} We claim that there is an action of $\LL_\Conf^+ M$ on $\ZZ_{U_P, \Pt^-}$ whose restriction to
    \[\Mod_M^{+, \vartheta} \mon \Gr_{M, X^\vartheta}\]
matches the usual action of $\LL_{X^\vartheta}^+ M$. This is because $\ZZ_{U_P, \Pt^-}$ admits a local description similar to that of $\ZZ_{N, \Pt^-}$. Namely, it is the base change of
    \[\Maps^\gen\left(\Disc, \pt \mon (U_P^- \bs G \git U_P) /M \right)\]
along the tautological section (\ref{eqn:tautological_section}). The left action of $M$ on $U_P^- \bs G \git U_P$ induces the action we seek. 

\subsubsection{First reduction} This model matches the strata
    \[\Mod_M^{+, \vartheta' \le \vartheta} \mon \Mod_M^{+, \theta}
    \qquad
    \HH_M^{+, \vartheta' \le \vartheta} \mon \HH_M^{+,\theta}.\]
Therefore, it suffices to show the singular support of $j_{\ZZ,!} \cdot k$ at every point of $\Mod_M^{+,\theta}$ is contained in the conormals to the strata $\Mod_M^{+, \vartheta' \le \vartheta}$.

\subsubsection{Second reduction} By the factorization of the Zastava diagram, we only need to consider points lying over the main diagonal $X \mon X^\theta$. These are the points belong to strata indexed by $\vartheta' \le (\theta)$ where $(\theta) \vdash \theta$ is the trivial partition.

\subsubsection{Third reduction} Since the Zastava diagram is etale local with respect to the curve, we may take $X = \AA^1$. Then the Zastava diagram admits a translation action by $\Ga$ commuting with the $\LL_\Conf^+ M$ action constructed above. The result follows because each stratum $\Mod_M^{+, \vartheta' \le (\theta)}$ is an orbit for the combined action of $\LL_\Conf^+ M \times \Ga$.

\subsection{End of the proof}

\subsubsection{} Let $\UU_\Sing \mon \Bunt_{P^-}$ be the open subset where the singular support of $j_{P,!}^-$ is transverse to the codifferential of $\Bunt_{P^-} \to \Bun_G$. The singular support bound implies:

\begin{proposition}
Let $x$ be a geometric point of ${}_\theta \Bunt_{P^-}$. If the image of $x$ under the projection ${}_\theta \Bunt_{P^-} \to \Bun_{P^-}$ belongs to the smooth locus of $\Bun_{P^-} \to \Bun_G$, then $x$ belongs to $\UU_\Sing$.
\end{proposition}

\subsubsection{Proof} We have the following Cartesian diagram:
\[
\begin{tikzcd}[
  column sep={1.5em},
  row sep={1.5em}
]
  \HH_M^{+, \vartheta' \le \vartheta} \arrow[r, hook] \arrow[dr, dashed] & \HH_M^{+, \theta} \arrow[d] & {}_\theta \Bunt_{P^-} \arrow[l] \arrow[d] & \\
  & \Bun_M & \Bun_{P^-} \arrow[l] \arrow[r] & \Bun_G
\end{tikzcd}
\]
Therefore it suffices to observe that diagonal dashed arrow is smooth.

\subsubsection{} Define
    \[\ZZt_\Sing = \UU_\Sing \underset{\Bunt_{P^-}}{\times} \ZZt.\]

\begin{proposition}
\label{proposition:Hecke_transitivity}
Every geometric point of $\ZZt$ can be taken to a point of $\ZZt_\Sing$ by the action of $\HHt$.
\end{proposition}

\subsubsection{Proof of Proposition \ref{proposition:Hecke_transitivity}} The proof comes in two parts.

\subsubsection{Part one} Let $z$ be a geometric point of $\ZZt$, and let $\FF_{P^-}$ be its image in $\Bun_{P^-}$. Then $\FF_{P^-}$ belongs to the smooth locus of $\Bun_{P^-} \to \Bun_G$ iff
\begin{equation}
\label{eqn:ULA_cohomology}  
    H^0\left(X, \FF_{P^-} \overset{P^-}{\times} \fg/\fp^- \right) \simeq 0
\end{equation}
To control this cohomology group, we choose an arbitrary reduction $\FF_{B^-}$ of $\FF_{P^-}$ to $B^-$ that is generically transverse to the image of $z$ in $\Bunc_N$. Then the vector bundle appearing in (\ref{eqn:ULA_cohomology}) acquires a complete flag whose graded pieces are the line bundles $\FF_T \overset{T}{\times} \alpha$ for $\alpha$ ranging over the roots appearing in $\fu_P$.

\subsubsection{Part two} From all of this data, we may extract the following defect divisors:
\begin{itemize}
    \item The image of $z$ in $\Conf_{G,P}$ is a divisor $D_1$, away from which $\FF_M$ has a degenerate $N_M$-reduction.
    \item The degeneracy of this $N_M$-reduction defines a $\Lambda_M$-valued divisor $D_2$ on $X \setminus |D_1|$.
    \item From the choice of $\FF_{B^-}$ we obtain a point of $\ZZ_{\overline{N},B^-}$ and hence a divisor $D_3$ in $\Conf_{G,B}$. 
\end{itemize}
Now consider a point $E$ of $\Conf_{G,P}$ whose support is disjoint from all of three of the above divisors, and let $E'$ be a lift of $E$ to $\Conf_{G,B}$. Perform a Hecke modification coming from
    \[
    \left(\ZZ_{N,B^-,\Conf}\right)_{E'} \to \left(\ZZ_{N,P^-,\Conf}\right)_E\]
to obtain a new point of $\ZZt$, say $z'$. The image of $z'$ in $\Bun_{P^-}$ comes with a reduction $\FF_{B^-}'$ to $B^-$, and we have
  \[\FF_T' \overset{T}{\times} \alpha \simeq \left(\FF_T \overset{T}{\times} \alpha \right)(\alpha(E'))\]
But we may choose $E$ so that the divisors $\alpha(E') = \alpha(E)$ for $\alpha$ appearing in $\fu_P$ are arbitrarily negative, thereby forcing the vanishing (\ref{eqn:ULA_cohomology}) for $\FF_{P^-}'$.

%% file: coJacquet.tex
\section{The co-Jacquet functor}
\label{section:co_Jac}

\subsubsection{} In this section we shall construct the following objects:
\begin{itemize}
    \item A Jacquet functor
        \[\Jac_{\co,\Conf} : \Whit(M)_\Conf \to \Whit(G)\]
    along with a commutative diagram
        \begin{equation}
        \label{eqn:coeff_Eis_Conf}
        \begin{tikzcd}[
          column sep={8em,between origins},
          row sep={4em,between origins}
        ]
          \DMod(\Bun_M) \arrow[r, "\Eis_!^-"] \arrow[d, "\coeff_*" left] & \DMod(\Bun_G) \arrow[d, "\coeff_*"] \\
          \Whit(M)_\Conf \arrow[r, "\Jac_{\co}"] & \Whit(G)
        \end{tikzcd}
        \end{equation}
    \item A duality
        \[\Whit(M)_\Conf^\vee \simeq \Whit(M)_\Conf\]
    along with a commutative diagram
        \begin{equation}
        \label{eqn:Jac_duality_Conf}
        \begin{tikzcd}[
          column sep={8em,between origins},
          row sep={4em,between origins}
        ]
          \Whit(G) \arrow[r, "\Jac"] \arrow[d, "\DD{[2 \delta_G]}" left] & \Whit(M)_\Conf \arrow[d, "\DD{[2\delta_M]}"] \\
          \Whit(G)^\vee \arrow[r, "\Jac_{\co}^\vee"] & \Whit(M)_\Conf^\vee 
        \end{tikzcd}
        \end{equation}
\end{itemize}

\subsection{Definition of the functor}

\subsubsection{} Let us return to the diagram (\ref{eqn:zastava_diagram}). We have, in particular, the following morphisms:
\[
\begin{tikzcd}[
    column sep={6em,between origins},
    row sep={2.5em,between origins}
    ]
    & \ZZt_{N,P^-} \arrow[dr, "u"] \arrow[dd, "w"] \arrow[dl, "v" above left] & \\
    \Bunc_N & & \Bunc_{N_M,\Conf}^\pol \arrow[dd, "\pi_\Conf"] \\
    & \Bunt_{P^-} & \\
    & & \Bun_M 
\end{tikzcd}
\]
Recall that
  \[\Eis_!^- = \qt_! \cdot \left( \widetilde{p}^*(-) \stimes j_{P,!}^-\right) \simeq \qt_* \cdot \left( \widetilde{p}^!(-) \shimes j_{P,!}^- \right) [2 \delta_M].\]
Therefore, we obtain a morphism
    \begin{align*}
    \coeff_* \cdot \Eis_!^- 
    & \simeq \Av_*^\Whit \cdot \aph\widetilde{p}_* \cdot \left( \aph\pi^! \cdot \qt^!(-) \shimes \aph\pi^! \cdot j_{P,!}^- \right) \, [2 \delta_M]\\
    & \qquad \to \Av_*^\Whit \cdot v_* \cdot \left( u^! \cdot \pi_\Conf^! (-) \shimes w^! \cdot j_{P,!}^- \right) \, [2 \delta_M]
    \end{align*}
The same analysis as Proposition \ref{proposition:zast_vanishing} shows that this map is an isomorphism. Therefore, we consider the diagram
\begin{equation}
\label{eqn:Jac_co_tilde}
\begin{tikzcd}[
          column sep={8.5em,between origins},
          row sep={4em,between origins}
        ]
  \DMod\left( \Bunc_{N_M, \Conf}^\pol \right) \arrow[r, "\Jact_\co"] \arrow[d, "\Av_*"] & \DMod(\Bunc_N) \arrow[d, "\Av_*"] \\
  \Whit(M)_\Conf \arrow[r, dashed, "\Jac_\co"] & \Whit(G)
\end{tikzcd}
\end{equation}
where
  \[\Jact_\co = v_* \cdot \left( u^!(-) \shimes w^! \cdot j_{P,!}^- \right) \, [2\delta_M].\]

\begin{proposition}
\label{proposition:Jac_co_factors}
The factorization $\Jac_\co$ indicated in the diagram (\ref{eqn:Jac_co_tilde}) exists.
\end{proposition}

\subsubsection{} Note that the functor $\Jac_\co$ is uniquely determined by the commutativity of (\ref{eqn:Jac_co_tilde}). This completes the construction of $\Jac_\co$ and the commutative diagram (\ref{eqn:coeff_Eis_Conf}) modulo the proof of Proposition \ref{proposition:Jac_co_factors}, which will appear in the next subsection.

\subsection{Duality} 

\subsubsection{} For $\WW \in \Whit(M)_\Conf$ the formula
  \[\langle \WW, - \rangle_M = \Gamma\left(\Bunc_{N_M,\Conf}^\pol, \WW \shimes - \right)\]
defines a continuous functor
  \[\DMod\left( \Bunc_{N_M,\Conf}^\pol \right) \to \Vect.\]
Furthermore, the map
    \[
    \left\langle \WW, \Av_*^{-\psi} \cdot \FF \right\rangle 
    \to \langle \WW, \FF \rangle_M 
    \]
is an isomorphism, and we obtain a perfect pairing
  \[\Whit_\psi(M)_\Conf \otimes \Whit_{-\psi}(M)_\Conf \to \Vect.\]

\begin{proposition}
For any $\WW_G \in \Whit(G)$ and $\FF_M \in \DMod\left( \Bunc_{N_M,\Conf}^\pol \right)$ we have
    \[
    \left\langle \Jac \cdot \WW_G, \FF_M \right\rangle_M 
    \simeq \left\langle \WW_G, \widetilde{\Jac}_{\co,\Conf} \cdot \FF_M \right\rangle_G [2 \delta_G - 2 \delta_M].
    \]
\end{proposition}

\subsubsection{Proof} This is clear because
\begin{align*}
  \LHS & \simeq \Gamma\left( \Bunc_{N(M),\Conf}^\pol, u_* \left\{ v^! \cdot \WW_G \shimes w^! \cdot j_{P,!}^- \right\} \shimes \FF_M \right) [2\delta_G] \\
  & \qquad \simeq \Gamma\left( \ZZt_{N,P^-}, v^! \cdot \WW_G \shimes w^! \cdot j_{P,!}^- \shimes u^! \cdot \FF_M \right) [2\delta_G] \\
  & \qquad\quad \simeq \Gamma\left( \Bunc_N, \WW_G \shimes v_* \left\{ w^! \cdot j_{P,!}^- \shimes u^! \cdot \FF_M \right\} \right) [2 \delta_G] \simeq \RHS.
\end{align*}

\subsubsection{Proof of Proposition \ref{proposition:Jac_co_factors}} Since the duality pairing factors through $\Av_*^{-\psi}$, the above proposition gives
  \[
  \left\langle \WW_G, \Av_*^{-\psi} \cdot \Jact_\co \cdot \FF_M \right\rangle_G
  \simeq \left\langle \Jac \cdot \WW_G, \Av_*^{-\psi} \cdot \FF_M \right\rangle_M [2\delta_M - 2\delta_G].
  \]
The result follows from the fact that $\langle -, - \rangle_G$ is a perfect pairing.

\subsubsection{Construction of (\ref{eqn:Jac_duality_Conf})} We have
    \[
    \langle \Jac \cdot \WW_G, \WW_M \rangle_M \, [2\delta_M]
    \simeq \langle  \WW_G, \Jac_\co \cdot \WW_M \rangle_G \, [2\delta_G]
    \]
for every $\WW_G \in \Whit(G)$ and $\WW_M \in \Whit(M)_\Conf$.

%% file: Ran.tex
\section{The factorizable Whittaker category}
\label{section:Ran}

\subsubsection{} In the section we will review the definitions of the Ran version $\Whit(G)_\Ran$ of the Whittaker category from \cite{loc_glob_whit} and relate it to the configuration space version $\Whit(G)_\Conf$. Although $\Whit(G)_\Ran$ is a naturally a unital factorization category, we will not make explicit use of this additional structure.

\subsection{Review of Whittaker categories}

\subsubsection{Local model} For a finite set $I$, we write
    \[\Whit(G)_I = \DMod(\Gr_{G,I})^{\LL N,\psi}  \qquad  \Whit(G)_{I,\co} = \DMod(\Gr_{G, I})_{\LL N, \psi}\]
By \cite[Appendix B]{lin_nearby_cycles}, these assignments define factorization categories which we denote
  \[\Whit(G)  \qquad  \Whit(G)_\co.\]
Furthermore, we have a canonical $\DMod(X^I)$-linear equivalence
  \[\Whit_\psi(G)_I^\vee \simeq \Whit_{-\psi}(G)_{I,\co}\]
inducing an equivalence
  \[\Whit_\psi(G)^\vee \simeq \Whit_{-\psi}(G)_\co\]
of factorization categories.

\subsubsection{Local duality} The functor of renormalized $*$-averaging
    \[\Av_{*,\ren} : \DMod(\Gr_{G,I}) \to \Whit(\Gr_{G,I})\]
induces an equivalence
    \[\PsId_\Whit : \Whit(G)_{I,\co} \simto \Whit(G)_I.\]
We obtain an equivalence of factorization categories
    \[\Whit_{-\psi}(G)^\vee \simto \Whit_\psi(G).\]

\subsubsection{Global model} Write $\Bunc_{N,I}^\pol$ for the stack classifying a point $x \in X^I$, and a principal $G$-bundle on $X$ equipped with a degenerate $N^\omega$-reduction away from $x$. We set
    \[\Whit(G)_I^\glob = \Whit\left(\Bunc_{N,I}^\pol \right).\]
This defines a sheaf of categories over the Ran space, which we denote
    \[\Whit(G)^\glob.\]
Note that $\Whit(G)^\glob$ acquires a weak factorization structure from the weak factorization structure
  \[\Bunc_{N,\Ran}^\pol \underset{\Ran}{\times} \left[ \Ran \times \Ran \right]_\disj \to \left[ \Bunc_{N,\Ran}^\pol \times \Bunc_{N,\Ran}^\pol \right]_\disj.\]
The comparison with the local category will show that $\Whit(G)^\glob$ is a genuine factorization category.

\subsubsection{Global duality} The perfect pairing
    \[
    \Whit_{-\psi}\left( \Bunc_{N,I}^\pol \right) \otimes \Whit_{\psi}\left( \Bunc_{N,I}^\pol \right) \to \Vect
    \qquad \langle \WW_1, \WW_2 \rangle = \Gamma\left( \Bunc_{N,I}^\pol, \WW_1 \shimes \WW_2 \right)\]
defines a factorizable equivalence
  \[\DD : \Whit_{\psi}(G)^\glob \simto \Whit_{-\psi}(G)^{\glob,\vee}.\]

\subsubsection{Local v global} The map
    \[\pi : \Gr_{G,I} \to \Bunc_{N,I}^{\pol}\]
induces equivalences
    \[
        \pi^! : \Whit(G)_I^\glob \to \Whit(G)_I  \qquad
        \pi_* : \Whit(G)_{\co,I} \to \Whit(G)_{\co, I}^\glob.
    \]
We use the same notation for the Ran versions. The interaction of this comparison with duality is summarized by the diagram below.
\begin{equation}
\label{eqn:Whit_duality_loc_global}
\begin{tikzcd}[
      column sep={5em,between origins},
      row sep={6em,between origins}
    ]
    \Whit_{-\psi}(G)^{\glob,\vee} & & & \Whit_\psi(G)^\glob \arrow[lll, "\DD" above] \arrow[dr, "\pi^!"] & \\
    \Whit_{-\psi}(G)^\vee \arrow[rr, "\sim"] \arrow[u, "(\pi^!)^\vee"] & & \Whit_\psi(G)_\co \arrow[rr, "\PsId{[-2 \delta_N]}"] \arrow[ur, "\pi_*"] & & \Whit_\psi(G)
\end{tikzcd}
\end{equation}

\subsubsection{Hecke action} 
\begin{itemize}
    \item We write $\Bun_{G,I}$ for $\Bun_G \times X^I$ and $\DMod(\Bun_G)_I$ for $\DMod(\Bun_{G,I})$, viewed as a $\DMod(X^I)$-module.
    \item Recall the Hecke stack
        \[
        \begin{tikzcd}[
      column sep={5em,between origins},
      row sep={2em,between origins}
    ]
          & \HH_{G,I} \arrow[dl, "\hl" above left] \arrow[dr, "\hr"] \arrow[dd, "s"] &    \\
          \Bun_G & & \Bun_G \\
          & X^I &
        \end{tikzcd}
        \]
    Note that $\DMod(\HH_{G,I})$ is naturally a $\DMod(X^I)$-algebra acting on $\DMod(\Bun_G)_I$.
    \item We write $\Sph_{G,I}$ for the spherical Hecke category over $X^I$. It is a $\DMod(X^I)$-algebra acting on $\DMod(\Bun_G)_I$ through the canonical monoidal functor
        \[
        \Sph_{G,I} \to \DMod(\HH_{G,I})   
        \qquad \SS \mapsto \HH(\SS).
        \]
    \item We write $\Sph_{G,\Ran}$ for the Ran version of $\Sph_{G,I}$. It is a monoidal factorization category.
\end{itemize}
Now $\Sph_{G,I}$ acts on all of the previously defined Whittaker categories over $X^I$. Let us discuss the compatibility of (\ref{eqn:Whit_duality_loc_global}) with the $\Sph_G$ action. The right triangle of (\ref{eqn:Whit_duality_loc_global}) is a diagram of right $\Sph_G$-modules, and the left arrow $(\pi^!)^\vee$ is a morphism of left $\Sph_G$-modules. The global Verdier duality is $\Sph_G$-linear with respect to the canonical anti-involution $\tau : \Sph_G \to \Sph_G^\rev$, as is the identification $\Whit_{-\psi}(G)^\vee \simto \Whit_\psi(G)_\co$.

\subsubsection{Whittaker coefficients and Poincare series} The functor of Whittaker coefficients
    \[
    \coeff_{I,*} : \DMod(\Bun_G)_I \overset{\pi^!}{\longrightarrow} \DMod\left(\Bunc_{N,I}^\pol \right) \overset{\Av_*}{\longrightarrow} \Whit(G)_I^\glob
    \]
is $\Sph_{G,I}$-linear. Therefore, its left adjoint
    \[\Poinc_{I,!} : \Whit(G)_I^\glob \to \DMod(\Bun_{G,I})\]
is automatically co-lax $\Sph_{G,I}$-linear. We write $\coeff_{\Ran,*}$ and $\Poinc_{\Ran,!}$ for the Ran variants of these functors.

\begin{proposition}
\label{proposition:Poinc_linear}
The functor $\Poinc_{I,!}$ is actually $\Sph_{G,I}$-linear.
\end{proposition}

\subsubsection{A triviality} We will use the following observation in both the proof of Proposition \ref{proposition:Poinc_linear} and in later sections without comment. Suppose that $f : M_1 \to M_2$ is a colax morphism of modules over a monoidal category $A$. Let us say that an object $m$ of $M_1$ is \emph{strict} if the colax constraint
  \[f(a \otimes m) \to a \otimes f(m)\]
is an isomorphism for all $a \in A$. Then the full subcategory of strict objects is a monoidal ideal.

\subsubsection{Proof of Proposition \ref{proposition:Poinc_linear}} Since $\Whit(G)_I$ is generated as a $\Sph_{G,I}$-module by the vacuum object
    \[\WW_{\vac, I} = \WW_\vac \boxtimes \omega_{X^I} \in \Whit\left(\Bunc_N \times X^I\right) \mon \Whit(G)_I,\]
it suffices to show that the colax constraint
    \[\Poinc_{I,!} \cdot (\SS \star \WW_{\vac,I})
    \to \SS \star \left( \Poinc_{I,!} \cdot \WW_{\vac,I} \right)\]
is an isomorphism. Now consider the following diagram.
\[
\begin{tikzcd}[
      column sep={4em,between origins},
      row sep={4em,between origins}
    ]
    & \Bunc_N \!\underset{\Bun_G}{\times}\! \HH_{G,I} \arrow[dl, "\hl" above left] \arrow[dr] \arrow[drrr, bend left=15, dashed, "\hr_{\!N}"] & & & \\
  \Bunc_N \arrow[dr, "\pi"] & & \HH_{G,I} \arrow[dl, "\hl"] \arrow[dr, "\hr \times s"] & & \Bunc_{N,I}^\pol \arrow[dl, "\pi_I"] \\
    & \Bun_G & & \Bun_{G,I} & 
\end{tikzcd}
\]
The colax constraint identifies with
\[
    \pi_{I,!} \cdot \hr_{\!N,!} \cdot \left( \WW_\vac \shoximes \SS_\HH \right)
    \to (\hr \times s)_* \cdot \left( \hl^! \cdot \pi_! \cdot \WW_\vac \shimes \SS_\HH \right),
\]
which is an isomorphism because $\hr \times s$ and $\hr_{\!N}$ are proper and $\SS_\HH$ is ULA over $\Bun_G$.

\subsubsection{} We write
    \[\Gamma(\Ran, -) : \DMod(\Bun_G)_\Ran \to \DMod(\Bun_G)\]
for pushforward along the projection
    \[\pr : \Bun_G \times \Ran \to \Bun_G.\]
It is left adjoint to $\pr^!$, and we have a commutative diagram
\[
\begin{tikzcd}[
      column sep={9em,between origins},
      row sep={4em,between origins}
    ]
    \DMod(\Bun_G)_\Ran^\vee \arrow[r, "\PsId_!"] \arrow[d, "(\pr^!)^\vee" left] & \DMod(\Bun_G)_\Ran \arrow[d, "{\Gamma(\Ran,-)}"] \\
    \DMod(\Bun_G)^\vee \arrow[r, "\PsId_!"] & \DMod(\Bun_G)
\end{tikzcd}
\]

\subsection{Conf v Ran}

\subsubsection{} The purpose of this subsection is to construct a functor
    \[\sprd_* : \Whit(M)_\Conf \to \Whit(M)_\Ran\]
along with the following compatibilities:
\begin{itemize}
    \item We have
        \begin{equation}
        \label{eqn:Poinc_Ran_sprd}
        \Gamma(\Ran, \Poinc_{\Ran,!} \cdot \sprd_*) \simeq \Poinc_{\Conf, !}.
        \end{equation}
    \item Let $\sprd^!$ denote the right adjoint to $\sprd_*$. Then there is a commutative diagram:
    \begin{equation}
    \label{eqn:sprd_duality}
    \begin{tikzcd}[
      column sep={9em,between origins},
      row sep={4em,between origins}
    ]
      \Whit(M)_\Conf^\vee \arrow[r, "\sim"] \arrow[d, "(\sprd^!)^\vee"] & \Whit(M)_\Conf \arrow[d, "\sprd_*"] \\
      \Whit(M)_\Ran^\vee \arrow[r, "\sim"] & \Whit(M)_\Ran
    \end{tikzcd}
    \end{equation}
\end{itemize}

\subsubsection{} The device relating $\Conf$ to $\Ran$ is a Cartesian diagram
\begin{equation}
\label{eqn:Conf_Ran_vac}
\begin{tikzcd}[
  column sep={7.5em,between origins},
  row sep={2em,between origins}
]
    & \Bunc_{N_M, \Ran}^{\pol, \neg} \arrow[dl, "\aph\pi" above left] \arrow[dr, "\aph i^\neg"] \arrow[dd] & \\
  \Bunc_{N_M, \Conf}^{\pol} \arrow[dd] & & \Bunc_{N_M, \Ran}^{\pol} \arrow[dd] \\
    & \Gr_{M^\ab, \Ran}^{\omega,\neg} \arrow[dl, "\pi" above left] \arrow[dr, "i^\neg"] &  \\
  \Conf_{G,P} & & \Gr^\omega_{M^\ab, \Ran}
\end{tikzcd}
\end{equation}
whose terms are defined as follows:
\begin{itemize}
    \item We write $\Gr_{M^\ab,\Ran}^\omega$ for the $\omega^\rhov$-twisted form of the affine Grassmannian. It classifies a point $x$ of $\Ran$, a principal $M^\ab$-bundle $\FF_{M^\ab}$, and a system of compatible maps
        \[
        \FFo_{M^\ab} \overset{M^\ab}{\times} \theta \to \omega^{\theta(\rhov)}
        \qquad
        \theta \in \Lambda_{G,P}.
        \]
    We write $\Gr_{M^\ab}^{\omega,\neg}$ for the closed sub-prestack where these maps extend to regular maps for all $\theta \in \Lambda_{G,P}^\pos$. When this happens we obtain a divisor $D \in \Conf_{G,P}$ supported on $x$. This defines the leftward map denoted $\pi$.
    \item The prestack $\Bunc_{N_M, \Ran}^\pol$ is the Ran version of $\Bunc_{N_M,I}^\pol$. It classifies a point $x$ of the Ran space, a principal $M$-bundle $\FF_M$ on $X$, and an $N_M^\omega$-Pl\"ucker data away from $x$. The right vertical map is defined by observing that such a Pl\"ucker data provides an identification between the bundles
        \[
        \FF_M \overset{M}{\times} M^\ab
        \qquad  \text{and} \qquad
        \omega^\rhov \overset{T}{\times} M^\ab
        \]
    away from $x$.
    \item We may define $\Bunc_{N_M, \Ran}^{\pol, \neg}$ as the fiber product of either face; the two are canonically identified.
\end{itemize}
Then we define
  \[\sprd_* = (\aph i^\neg)_* \cdot \aph\pi^!.\]

\subsubsection{} 
\label{subsection:sprd_compatibilities}
In fact, the map
    \[\Gr_{M^\ab}^{\omega,\neg} \mon \Conf \times \Ran\]
identifies $\Gr_{M^\ab}^{\omega,\neg}$ with the locus of pairs $(D, x)$ cut out by the closed condition $|D| \mon x$. Let us use this to obtain the desired compatibilities.
\begin{itemize}
    \item Since $\pi$ is pseudo-ind-proper and $i^\neg$ is a closed immersion, $\sprd_*$ admits a continuous right adjoint
        \[\sprd^! \simeq \aph\pi_* \cdot (\aph i^\neg)^!.\]
    This gives the compatibility (\ref{eqn:sprd_duality}) with duality.
    \item The compatibility (\ref{eqn:Poinc_Ran_sprd}) with Poincare series follows from the universal homological contractibility of $\pi$.
\end{itemize}

%% file: Ran_upgrades.tex
\section{Ran upgrades}
\label{section:Ran_upgrades}

\subsubsection{} In this section we will upgrade our previous constructions to the Ran setting.

\subsection{The Ran version of $\eta_G$}

\subsubsection{} For convenience, we will write
  \[\Poinc_? = \PsId_! \cdot \coeff_*^\vee \cdot \PsId_\Whit^{-1}.\]
We write $\Poinc_{\Ran,?}$ and $\Poinc_{I,?}$ the variants of this functor over $\Ran$ and $X^I$.

\subsubsection{} Consider the diagram
\[
\begin{tikzcd}[
  column sep={8em,between origins},
  row sep={4em,between origins}
]
  \Whit(G)_\Ran^\vee \arrow[r] \arrow[d, "\coeff_*^\vee" left] & \Whit(G)_\Ran \arrow[d, "\Poinc_!"] \\
  \DMod(\Bun_G)^\vee \arrow[r, "\PsId_!"] & \DMod(\Bun_G)
\end{tikzcd}
\]
of $\Sph_{G,\Ran}$-linear categories. The goal of this section is to extend $\eta_G$ to a map
  \[\eta_{G,\Ran} : \Poinc_{\Ran,!} \to \Poinc_{\Ran, ?} [2\delta_G]\]
of $\Sph_{G,\Ran}$-linear functors. We will write the construction over each power of the curve.

\subsubsection{} In fact we can simply repeat the definition of $\eta_G$ in order to define a map
  \[\eta_{G,I} : \Poinc_{I,!} \to \Poinc_{I,?} [2\delta_G]\]
of plain functors. That is, we consider the diagram
\[
\begin{tikzcd}[
  column sep={12em,between origins},
  row sep={4em,between origins}
]
  \left( \Bunc_{N,I}^\pol \underset{\Bun_G}{\times} \Bunc_G \right)^\trans \arrow[r, hook] & \Bunc_{N,I}^\pol \underset{\Bun_G}{\times} \Bunc_G \arrow[d, "\pi_2"] \\
   & \Bun_G
\end{tikzcd}
\]
and define $\eta_{G,I}$ to be
\[
\begin{tikzcd}[
  column sep={6em,between origins},
  row sep={4em,between origins}
]
  & \pi_{2, !} \cdot \left( \WW \stoximes r_{G,!} \right)^\trans \arrow[dr] \arrow[dl, "\sim"] & & \pi_{2, *} \cdot \left( \WW \shoximes r_{G,!} \right) \arrow[dl, "\sim"] \\
  \pi_{2, !} \cdot \left( \WW \stoximes r_{G,!} \right) [2\delta_G] & & \pi_{2, *} \cdot \left( \WW \shoximes r_{G,!} \right)^\trans [2\delta_G]&
\end{tikzcd}
\]
The left arrow is obviously an isomorphism, and the right arrow is an isomorphism by the same argument as Proposition \ref{proposition:cis_vanishing}.

\begin{proposition}
\label{proposition:Sph_linear}
The map $\eta_{G,I}$ is canonically $\Sph_{G,I}$-linear.
\end{proposition}

\subsubsection{Remarks}
\begin{itemize}
    \item[(i)] This result is not tautological. It requires the functorial construction of commutative diagrams:
\[
\begin{tikzcd}[
  column sep={10em,between origins},
  row sep={4em,between origins}
]
  \SS \star ( \Poinc_! \cdot \WW) \arrow[r] \arrow[d, "\sim"] & \SS \star (\Poinc_? \cdot \WW) \arrow[d, "\sim"] \\
  \Poinc_! \cdot (\SS \star \WW) \arrow[r] & \Poinc_? \cdot (\SS \star \WW)
\end{tikzcd}
\]
The right vertical arrow uses the fact that $\PsId_{G,I,!}$ is linear for the action of $\Sph_{G,I}$. It is easy to prove this by observing that 
\[
    \PsId_{G,I,!} \star \SS \simeq \SS \star \PsId_{G,I,!} \in \DMod(\Bun_G \times \Bun_G \times X^I)
\]
is given by the $!$-pushforward of $\SS$ along the map
  \[\hl \times \hr \times s : \HH_{G,I} \to \Bun_G \times \Bun_G \times X^I.\]
But this identification is difficult to use because has the wrong (i.e. $*$-pull and $!$-push) nature. Therefore, we should remedy this by constructing a compactification of $\hl \times \hr \times s$ as an $(\HH_{G, I}, \HH_{G,I})$-bimodule.
    \item[(ii)] Since $\Whit(G)_\Ran$ is freely generated by $\WW_\vac$ under the action of $\Rep(\LG)_\Ran$, we automatically obtain a $\Rep(\LG)_\Ran$-linear version of $\eta_{G,\Ran}$. But this map is not obviously $\Sph_{G, \Ran}$-linear. Furthermore, it is not clear that this na\"ively defined map agrees with $\eta_{G,\Ran}$ as a plain natural transformation and this identification will be needed in the sequel.
\end{itemize}

\subsubsection{} We write
  \[\Vin\Bun_{G,I}^\pol \to \Bun_G \times \Bun_G \times X^I\]
for the stack classifying a point $x \in X^I$ along with the following mapping data:
\[
\begin{tikzcd}[
  column sep={8em,between origins},
  row sep={4em,between origins}
]
  \eta \arrow[r, hook] \arrow[d] & X \setminus x \arrow[r, hook] \arrow[d] & X \arrow[d] \\ 
  G \bs \Vin^\circ_G / G \arrow[r, hook] & G \bs \Vin_G / G \arrow[r] & G \bs T_\ad^+ / G
\end{tikzcd}
\]
We remind the reader that the action of $G \times G$ on $T_\ad^+$ appearing the bottom right of this diagram is trivial. As with $\Vin\Bun_G$, this stack comes with an action of $T$ and we define
  \[\Bunc_{G,I}^\pol = \Vin\Bun_{G,I}^\pol / T.\]
Here are the salient geometric properties of this stack:

\begin{proposition}
\label{proposition:Hecke_comp_propoerties}
We have the following:
\begin{itemize}
    \item[(i)] There is a factorization
      \[\Bunc_{G,I} \mon \Bunc_{G,I}^\pol \to \Bun_G \times \Bun_G \times X^I\]
    where the first map is a closed embedding and the second is ind-proper.
    \item[(ii)] There is an open embedding
      \[r_{G,I}^\pol : \HH_{G,I} \times \BB Z_G \mon \Bunc_{G,I}^\pol.\]
    \item[(iii)] The Hecke stack
      \[\HH_{G,I} \underset{X^I}{\times} \HH_{G,I} \to \Bun_G \times \Bun_G \times X^I\] 
    acts on $\Bunc_{G,I}^\pol$, and both action maps
      \[\HH_{G,I} \underset{\Bun_{G,I}}{\times} \Bunc_{G,I}^\pol \to \Bunc_{G,I}^\pol  \qquad  \Bunc_{G,I}^\pol \underset{\Bun_{G,I}}{\times} \HH_{G,I} \to \Bunc_{G,I}^\pol\]
    are ind-proper.
\end{itemize}
\end{proposition}

\subsubsection{} The verification of Proposition \ref{proposition:Hecke_comp_propoerties} is straightforward, but we will include the details in the next subsection. For now let us simply describe the fiber of $\Vin\Bun_{G,I}^\pol$ over a point $x \in X^I$ in the case $G = \SL_2$. It classifies a pair of principal $\SL_2$-bundles along with a nonzero rational map
    \[\alpha : \EE_1 \dashrightarrow \EE_2\]
of coherent sheaves which is regular away from $x$. 

\begin{proposition}
The functor
\[
    \Sph_{G,I} \to \DMod\left(\Bunc_{G,I}^\pol \right)
    \qquad \SS \mapsto r_{G,I,!}^\pol \cdot \HH(\SS)
\]
is canonically a map of $\Sph_{G,I}$-bimodules over $X^I$.
\end{proposition}

\subsubsection{} In what follows, we will mainly be interested in the full subcategory 
    \[\DMod\left(\Bunc_{G,I}^\pol \right)^\good \mon \DMod\left(\Bunc_{G,I}^\pol \right)\]
generated by the image of this functor.

\subsubsection{Proof} The functor in question is automatically colax monoidal, so it suffices to show that the colax morphisms
  \[r_{G,I,!}^\pol(\SS \star \id_{\Sph_{G,I}}) \to \SS \star r_{G,I,!}^\pol(\id_{\Sph_{G,I}})  \qquad  r_{G,I,!}^\pol(\id_{\Sph_{G,I}} \star \SS) \to r_{G,I,!}^\pol(\id_{\Sph_{G,I}}) \star \SS\]
are isomorphisms. Indeed, this would imply that every object of $\Sph_{G,I}$ is strict for both the left and right actions. This then implies that the functor is strictly linear.

\subsubsection{} By symmetry, we only need to consider the right action. Consider the diagram
\[
\begin{tikzcd}[
  column sep={4em,between origins},
  row sep={4em,between origins}
]
    & \HH_{G,I} \arrow[dl] \arrow[rr, "\aph r_G"] & & \Bunc_G \!\underset{\Bun_G}{\times}\! \HH_{G,I} \arrow[dl] \arrow[dr] & \\
    \Bun_G \arrow[rr, "r_G"] & & \Bunc_G \arrow[dr] & & \HH_{G,I} \arrow[dl] \\
     & & & \Bun_G &
\end{tikzcd}
\]
The colax constraint for the right action is the pushforward of
  \[{\aph r_G}_! \cdot \HH(\SS) \to r_{G,!} \shoximes \HH(\SS) \]
along the proper morphism
\[
    \Bunc_G \underset{\Bun_G}{\times} \HH_{G,I} 
    \simeq \Bunc_{G,I} \underset{\Bun_{G,I}}{\times} \HH_{G,I}
    \to \Bunc_{G,I}^\pol.
\]
This map is an isomorphism because $\HH(\SS)$ is ULA over $\Bun_G$.

\subsubsection{Notation} 
\begin{itemize}
    \item We write $(-)_\HH$ and $(-)^\trans$ for restriction to the open loci
\[
    \Bunc_{N,I}^\pol \underset{\Bun_{G,I}}{\times} \HH_{G,I}
    \mon \left(\Bunc_{N,I}^\pol \underset{\Bun_{G,I}}{\times} \Bunc_{G,I}^\pol \right)^\trans 
    \mon \Bunc_{N,I}^\pol \underset{\Bun_{G,I}}{\times} \Bunc_{G,I}^\pol.
\]
    \item We have the ind-representable morphism
    \[
    \pi_2 : \left(\Bunc_{N,I}^\pol \underset{\Bun_{G,I}}{\times} \Bunc_{G,I}^\pol \right)^\trans
    \to \Bun_{G,I}.
    \]
\end{itemize}

\subsubsection{} Now define a pair of functors
\[
    \Whit(G)_I^\glob \otimes \DMod\left(\Bunc_{G,I}^\pol \right)^\good
    \to \DMod(\Bun_{G,I}).
\]
by
\[
    \Poinc_{I,!}^\ext(\WW, \FF) = \pi_{2,!} \left( \WW \shoximes \FF \right)_\HH \qquad
    \Poinc_{I,?}^\ext(\WW, \FF) = \pi_{2,*} \left( \WW \shoximes \FF \right)^\trans
\]
By forgetting supports, we obtain a natural transformation
  \[\eta_{G,I}^\ext : \Poinc_{I,!}^\ext \to \Poinc_{I,*}^\ext.\]

\begin{proposition}
We have:
\begin{itemize}
    \item[(i)] The functor $\Poinc_{I,!}^\ext$ descends to a $\Sph_{G,I}$-linear functor
        \[\Whit(G)_I^\glob \underset{\Sph_{G,I}}{\otimes} \DMod\left(\Bunc_{G,I}^\pol \right)^\good \to \DMod(\Bun_{G,I})\]
    \item[(ii)] The functor $\Poinc_{I,*}^\ext$ descends to a $\Sph_{G,I}$-linear functor
      \[\Whit(G)_I^\glob \underset{\Sph_{G,I}}{\otimes} \DMod\left(\Bunc_{G,I}^\pol \right)^\good \to \DMod(\Bun_{G,I}).\]
    \item[(iii)] Furthermore, $\eta_{G,I}^\ext$ descends to a map of $\Sph_{G,I}$-linear functors.
\end{itemize}
\end{proposition}

\subsubsection{Proof of (i)} Note that $\Poinc_{I,!}^\ext$ admits a factorization
\[
    \Whit(G)_I^\glob \otimes \DMod\left( \Bunc_{G,I}^\pol \right)^\good
    \to \Whit(G)_I^\glob \otimes \DMod(\HH_{G,I})^\good
    \to \DMod(\Bun_{G,I})
\]
where 
    \[\DMod(\HH_{G,I})^\good \mon \DMod(\HH_{G,I})\]
denotes the full subcategory generated by the image of
    \[\HH(-) : \Sph_{G,I} \to \DMod(\HH_{G,I}).\]
The first functor in this factorization is linear for all three $\Sph_{G,I}$ actions, and the second functor factors through a $\Sph_{G,I}$-linear functor
\[
    \Whit(G)_I^\glob \underset{\Sph_{G,I}}{\otimes} \DMod(\HH_{G,I})^\good
    \to \DMod(\Bun_{G,I}).
\]

\subsubsection{Proof of (ii)} Actually the analogous claim for all of $\DMod\left(\Bunc_{G,I}^\pol \right)$ holds by base change.

\subsubsection{Proof of (iii)} This follows from the compatibility of forgetting supports with proper pushforwards.

\subsubsection{Proof of Proposition \ref{proposition:Sph_linear}} Form the map
\[
    \Whit(G)_I^\glob \simeq \Whit(G)_I^\glob \underset{\Sph_{G,I}}{\otimes} \Sph_{G,I} 
    \to \Whit(G)_I^\glob \underset{\Sph_{G,I}}{\otimes} \DMod\left(\Bunc_{G,I}^\pol \right)^\good
\]
of right $\Sph_{G,I}$-modules. Its composite with $\eta_{G,I}^\ext : \Poinc_{I,!}^\ext \to \Poinc_{I,*}^\ext$ identifies with $\eta_{G,I}$.

\subsection{Proof of Proposition \ref{proposition:Hecke_comp_propoerties}}

\subsubsection{Proof of (i)} To simplify notation, we will write the proof for a fixed point $x \in X^I$. Let $\UU = X \setminus x$. We claim that that both of the maps
    \[
    \Vin\Bun_G \to \Vin\Bun_G^\pol \to \Vin\Bun_G(\UU)
    \]
are closed immersions. This follows from the fact that the morphism
    \[\Maps(X, S/H) \to \Maps(\UU, S/H)\]
is a closed immersion for any finite scheme $S$ with a group action by $H$. Indeed, the second arrow is a base change of
    \[\Maps(X, G \bs T_\ad^+ / G) \to \Maps(\UU, G \bs T_\ad^+ / G)\]
while the composite is a base change of
    \[\Maps(X, G \bs \Vin_G / G) \to \Maps(\UU, G \bs \Vin_G / G).\]
The first arrow is then a closed immersion by the two out of three property. Now take quotients by $T$ to obtain the closed immersions
    \[\Bunc_G \mon \Bunc_G^\pol \mon \Bunc_G(\UU).\]
Replacing a complete curve with an open curve in the proof \cite[Appendix A]{VinGr} that $\Bunc_G$ is proper over $\Bun_G \times \Bun_G$ shows that $\Bunc_G(\UU)$ is ind-proper over $\Bun_G(\UU) \times \Bun_G(\UU)$. Hence $\Bunc_G^\pol$ is proper over the closed substack
    \[\Bun_G \times \Bun_G \mon \Bun_G(\UU) \times \Bun_G(\UU).\]

\subsubsection{Proof of (ii)} Note that $\HH_{G,I}$ identifies with the base change of
    \[\Vin\Bun_{G,I}^\pol \to \Maps(X, T_\ad^+) \simeq T_\ad^+\]
along the open subset $T_\ad \mon T_\ad^+$. Taking the quotient by $T$ gives $r_{G,I}^\pol$.

\subsubsection{Proof of (iii)} The existence of the action in (iii) is clear and the ind-properness follows from that of the affine Grassmannian.

\subsection{The Ran version of the Jacquet functor}

\subsubsection{} The purpose of this subsection is to construct the Ran version of the Jacquet functor along with some of its attendant compatibilities. Recall that we have the space
    \[\Ran_\Ran \mon \Ran \times \Ran\]
classifying pairs of points $x \subset y$ in $\Ran \times \Ran$. The Ran version of the Jacquet functor will be a map        \[\Jac_\Ran : \Whit(G)_\Ran \to \Whit(M)_{\Ran_\Ran},\]
where $\Whit(M)_{\Ran_\Ran}$ is the pullback of $\Whit(M)$ along the projection onto the second copy of $\Ran$. 

\subsubsection{} We will work over powers of the curve. As before, the definition of this functor will actually be as a composite of functors
    \[\Jac_I : \Whit(G)_I \to \Whit(M)_{\Conf_I} \to \Whit(M)_{\Ran_I}\]
Here is a complete list of the objects we need to construct in this subsection:
\begin{itemize}
    \item A Jacquet functor
        \[\Jac_{\Conf,I} : \Whit(G)_I \to \Whit(M)_{\Conf_I}\]
    with the compatibility
        \begin{equation}
        \label{eqn:Jac_Conf_I_Poinc}
        \begin{tikzcd}[
          column sep={8em,between origins},
          row sep={4em,between origins}
        ]
            \Whit(G)_I \arrow[r, "\Jac"] \arrow[d, "\Poinc_!"] & \Whit(M)_{\Conf_I} \arrow[d, "\Poinc_!"] \\
            \DMod(\Bun_G)_I \arrow[r, "\CT_!^-"] & \DMod(\Bun_M)_I
        \end{tikzcd}
        \end{equation}
    \item A co-Jacquet functor    
        \[\co\Jac_{\Conf,I} : \Whit(M)_{\Conf_I} \to \Whit(G)_I\]
    with the compatibility
        \begin{equation}
        \label{eqn:Jac_co_Conf_I_coeff}
        \begin{tikzcd}[
          column sep={8em,between origins},
          row sep={4em,between origins}
        ]
            \DMod(\Bun_M)_I \arrow[r, "\Eis_!^-"] \arrow[d, "\coeff_*"] & \DMod(\Bun_G)_I \arrow[d, "\coeff_*"] \\
            \Whit(M)_{\Conf_I} \arrow[r, "\Jac_\co"] & \Whit(G)_I
        \end{tikzcd}
        \end{equation}        
    \item A duality
        \[\DD : \Whit_\psi(M)_{\Conf_I} \simto \Whit_{-\psi}(M)_{\Conf_I}^\vee\]
    with the compatibility
        \begin{equation}
        \label{eqn:Jac_I_duality}
        \begin{tikzcd}[
          column sep={10em,between origins},
          row sep={4em,between origins}
        ]
            \Whit_\psi(G)_I \arrow[r, "\Jac"] \arrow[d, "\DD{[2\delta_G]}" left] & \Whit_\psi(M)_{\Conf_I} \arrow[d, "\DD{[2\delta_M]}"] \\
            \Whit_{-\psi}(G)_I^\vee \arrow[r, "\Jac_\co^\vee"] & \Whit_{-\psi}(G)_{\Conf_I}^\vee
        \end{tikzcd}
        \end{equation}
    \item A functor
        \[\sprd_{I,*} : \Whit(M)_{\Conf_I} \to \Whit(M)_{\Ran_I}\]
    with the compatibility
        \begin{equation}
        \label{eqn:sprd_I_duality}
        \begin{tikzcd}[
          column sep={11em,between origins},
          row sep={4em,between origins}
        ]
            \Whit_\psi(M)_{\Conf_I} \arrow[r, "\sprd_*"] \arrow[d, "\DD" left] & \Whit_\psi(M)_{\Ran_I} \arrow[d, "\DD"] \\
            \Whit_{-\psi}(M)_{\Conf_I}^\vee \arrow[r, "(\sprd^!)^\vee"] & \Whit_{-\psi}(M)_{\Ran_I}^\vee
        \end{tikzcd}
        \end{equation}
    and an identification
        \begin{equation}
        \label{eqn:sprd_I_Poinc}
            \Poinc_{\Conf_I,!} \simeq \Poinc_{\Ran_I,!} \cdot \sprd_*
        \end{equation}
    of functors
        \[\Whit(M)_{\Conf_I} \to \Whit(M)_{\Ran_I}.\]
    Here $\sprd_I^!$ is the right adjoint to $\sprd_{I,*}$.
\end{itemize}

\subsubsection{Pole-adding} Let $\YY \to X^I$ be any prestack equipped with an action of the constant monoid scheme $(\Lambda^\pos)^I$. We will think of and refer to such actions as \emph{zero-adding}. Inverting this action gives a new prestack denoted $\YY^\pol$. All of our examples will have the feature that the action maps
  \[\lambda : \YY \to \YY  \qquad \lambda \in (\Lambda_{G,P}^\pos)^I.\]
are closed immersions.

\subsubsection{} Let us list a number of actions to which we will apply the pole-adding procedure. In many cases, the action will factor through the map $\Lambda^\pos \epi \Lambda_{G,P}^\pos$.
\begin{itemize}
    \item We give $\Bun_{G,I}$ and $\Bunt_{P^-, I}$ the trivial actions. So we will not write $\Bun_{G,I}^\pol$ or $\Bunt_{P^-, I}^\pol$.
    \item We give $\Conf_{G,P} \times X^I$ the $(\Lambda_{G,P}^\pos)^I$-action of
        \[\lambda \cdot (D, x) = (D + \lambda \cdot x, x).\]
    To reduce notation, we write
        \[\Conf_I = (\Conf_{G,P} \times X^I)^\pol.\]
    \item We give $\Bunc_{N,I}$ the usual action of $(\Lambda^\pos)^I$ by adding zeros. Note that the map $\Bunc_{N,I} \to \Bun_{G,I}$ is equivariant for this action. We warn that this action is not compatible with the formation of Whittaker categories.
    \item In the same way,
        \[\Bunc_{N_M, \Conf}^\pol \times X^I\]
    has an action of $\Lambda_{G,P}^\pos$ by adding zeros for which the projection to $\Conf \times X^I$ is equivariant. We write
        \[\Bunc_{N_M, \Conf_I}^\pol \to \Conf_I\]
    for the result of applying pole-adding to this map.
    \item The open substack
        \[
        \ZZ_{\Nc, \Pt^-, I}
        \mon \left(\Bunc_N \underset{\Bun_G}{\times} \Bunt_{P^-}\right)_I
        \]
    is stable under zero-adding, and the map
        \[
        \ZZ_{\Nc, \Pt^-, I}
        \to \Bunc_{N_M, \Conf}^\pol \times X^I
        \]
    is equivariant for this action.
\end{itemize}

\subsubsection{} The pole-adding procedure then gives a polar variant of (\ref{eqn:zastava_diagram}).
\[
\begin{tikzcd}[
  column sep={4em,between origins},
  row sep={4em,between origins}
]
     & \ZZ_{\Nc, \Pt^-, I}^\pol \arrow[dl, "v" above left] \arrow[dr, "t"] \arrow[drrr, dashed, bend left=15, "u"] & &  & \\
    \Bunc_{N,I}^\pol \arrow[dr] & & \Bunt_{P^-, I} \arrow[dl] \arrow[dr] & & \Bunc_{N_M, \Conf_I}^\pol \arrow[dl] \\
     & \Bun_{G,I} & & \Bun_{M,I} &
\end{tikzcd}
\]
We define
\[
    \Jac_\Conf : \Whit(G)_I \to \Whit(M)_{\Conf_I}
    \qquad \WW \mapsto u_* \cdot \left( \WW \shoximes j_{P,!}^- \right)_\ZZt [2 \delta_G]
\]
where the formation of $\shoximes$ is with respect to the fiber square
    \[
    \begin{tikzcd}[
  column sep={3em,between origins},
  row sep={3em,between origins}
]
        & \square \arrow[dl] \arrow[dr] \\
      \Bunc_{N,I} \arrow[dr] & & \Bunt_{P^-} \arrow[dl] \\
       & \Bun_G
    \end{tikzcd}
    \]
The fact that this formula produces an object in the Whittaker category is checked in the same way as Proposition \ref{proposition:Jac_Whit_preservation}. Now we have as before the natural transformations
\begin{align*}
    \CT_{I,!}^- \cdot \Poinc_{G,I,!}
    & \overset{\alpha}{\longleftarrow} \pi_! \cdot u_! \left( - \stoximes j_{P,!}^- \right)_\ZZt \\
    & \qquad\qquad \overset{\beta}{\longrightarrow} \pi_! \cdot u_* \left( - \stoximes j_{P,!}^- \right)_\ZZt 
    \overset{\gamma}{\longrightarrow} \Poinc_{M,I,!} \cdot \Jac_{\Conf,I}.
\end{align*}

\begin{proposition}
All of these maps are isomorphisms.
\end{proposition}

\subsubsection{Proof} The vanishing argument used to prove Proposition \ref{proposition:zast_vanishing} shows that $\alpha$ is an isomorphism. Also $\beta$ is an isomorphism because $u$ is proper: it was obtained from a proper map via the pole-adding procedure. We will prove that $\gamma$ is an isomorphism by use of the following reductions.

\subsubsection{First reduction: to a point} The four functors appearing above are colax linear over $\DMod(X^I)$. We claim that they are actually strictly linear. For $\CT_{I,!}^- \cdot \Poinc_{G,I,!}$, this follows from the second adjointness $\CT_{I,!}^- \simeq \CT_{I,*}$ and the fact that $\Poinc_{G,I,!}$ is even $\Sph_{G,I}$-linear. Since $\alpha$ and $\beta$ are isomorphisms, it remains to observe $\Poinc_{M,I,!} \cdot \Jac_\Conf$ is linear because both of its factors are. Therefore, to check that $\gamma$ is an equivalence, we may do so stratawise.

\subsubsection{Second reduction: to $\WW_{\vac, x}$} The restriction of $\Whit(G)_I$ to each stratum of $X^I$ is actually constant, so it suffices to check that $\gamma$ is an equivalence pointwise. For this it is enough to show that the purity morphism
    \[
    \left( - \stoximes j_{P,!}^- \right)_\ZZt 
    \to \left( - \shoximes j_{P,!}^- \right)_\ZZt \, [2 \delta_{G}]
    \]
defines an isomorphism of functors
    \[
    \Whit(G)_x \to \DMod\left(\ZZ_{\Nc,\Pt,x}^\pol \right).
    \]
The first functor is $\Sph_{G,x}$-linear by the ULA property of Hecke operators and the second is obviously so. Therefore, the result follows from Lemma \ref{lem:zast_ULA}, which says that the purity morphism is an isomorphism on the object $\WW_{\vac, x}$.

\subsubsection{Dualities} 

As before, we have a $\DMod(X^I)$-linear perfect pairing
\[
    \Whit_\psi(M)_{\Conf_I}^\glob \otimes \Whit_{-\psi}(M)_{\Conf_I}^\glob \to \DMod(X^I)
\]
given by $!$-tensoring the underlying sheaves followed by $*$-pushforward along the map
  \[\Bunc_{N, \Conf_I}^\pol \to \Conf_I \to X^I.\]
We take $\Jac_\co$ to be the descent
\[
\begin{tikzcd}[
  column sep={9em,between origins},
  row sep={4em,between origins}
]
    \DMod\left( \Bunc_{N,I}^\pol \right) \arrow[r] \arrow[d] & \DMod\left( \Bunc_{N_M, \Conf_I}^\pol \right) \arrow[d] \\
    \Whit(G)_I \arrow[r, dashed, "\Jac_\co"] & \Whit(M)_{\Conf_I} 
\end{tikzcd}
\]
where the upper arrow is
  \[\WW \mapsto v_* \left( t^! \cdot j_{P,!}^- \shimes u^! \cdot \WW \right)[2 \delta_{M}].\]
It is exchanged with $\Jac$ under the above duality.

\subsubsection{Conf v Ran} We need a version of the diagram (\ref{eqn:Conf_Ran_vac}) with marked points.
\[
\begin{tikzcd}[
  column sep={8em,between origins},
  row sep={2.5em,between origins}
]
    & \Bunc_{N_M, \Ran_I}^{\neg, \pol} \arrow[dl, "\aph\pi"] \arrow[dr, "\aph i^\neg"] \arrow[dd] & \\
  \Bunc_{N_M, \Conf_I}^{\pol} \arrow[dd] & & \Bunc_{N_M, \Ran_I}^{\pol} \arrow[dd] \\
    & \Gr_{M^\ab, \Ran_I}^{\omega,\neg, \pol} \arrow[dl] \arrow[dr] &  \\
  \Conf_I & & \Gr_{M^\ab, \Ran_I}^\omega
\end{tikzcd}
\]
This diagram is constructed by applying the pole-adding procedure to the diagram
\[
\begin{tikzcd}[
  column sep={8em,between origins},
  row sep={2.5em,between origins}
]
    & \square \arrow[dl] \arrow[dr] \arrow[dd] & \\
  \Bunc_{N_M, \Conf}^\pol \times X^I \arrow[dd] & & \Bunc_{N_M, \Ran_I}^{\pol} \arrow[dd] \\
    & \Gr_{M^\ab, \Ran_I}^{\omega,\neg} \arrow[dl, "\pi" above left] \arrow[dr] &  \\
  \Conf \times X^I & & \Gr^\omega_{M^\ab, \Ran_I}
\end{tikzcd}
\]
Here are the definitions:
\begin{itemize}
    \item The right face of this diagram is the base change of the right face of (\ref{eqn:Conf_Ran_vac}) along the map $\Ran_I \to \Ran$. The right vertical arrow carries an action of the group scheme $\Gr_{M^\ab, X^I}$, and we let $(\Lambda_{G,P}^\pos)^I$ act on this arrow through the map
      \[(\Lambda_{G,P}^\pos)^I \times X^I \to \Gr_{M^\ab, X^I}\]
    of monoids over $X^I$. This $(\Lambda_{G,P}^\pos)^I$-action preserves the closed sub-prestack $\Gr_{M^\ab, \Ran_I}^{\omega,\neg}$. This completes the definition of the right face and shows that the adding poles procedure does not affect the right vertical arrow.
    \item The leftward map $\pi$ is given by the projections
      \[
      \Gr_{M^\ab, \Ran_I}^{\omega,\neg} \to \Gr_{M^\ab, \Ran}^{\omega,\neg} \to \Conf
      \qquad
      \Gr_{M^\ab, \Ran_I}^\neg \to \Ran_I \to X^I.
      \]
    As before, the fiber products of the two faces are canonically identified.
\end{itemize}
Now we define
    \[
    \sprd_{I,*} : \Whit(M)_{\Conf_I} \to \Whit(M)_{\Ran_I}
    \qquad \WW \mapsto (\aph i^\neg)_* \cdot \aph\pi^! \cdot \WW.
    \]
One obtains the compatibilities (\ref{eqn:sprd_I_duality}) and (\ref{eqn:sprd_I_Poinc}) exactly as in \S\ref{subsection:sprd_compatibilities}.

%% file: 3_morphism.tex
\section{Proof of the main theorem} 
\label{section:3_morphism}

\subsubsection{} In this final section we will prove the main theorem in the following form:

\begin{theorem}
\label{thm:theorem_Ran}
The natural transformation $\eta_{G,\Ran}$ is an isomorphism.
\end{theorem}

\subsection{Proof of Theorem \ref{thm:theorem_Ran}}

\subsubsection{} Our proof proceeds by induction on $G$. We make a series of reductions.

\subsubsection{First reduction} We already know that $\eta_{G,\Ran}$ has Eisenstein cofiber. Therefore, it suffices to show that the cofiber of $\eta_{G, \Ran}$ is cuspidal. By Hecke linearity, it suffices to show that the cofiber of
    \[\eta_G \simeq \eta_{G,\Ran}(\WW_\vac)\]
is cuspidal.

\subsubsection{Second reduction} Fix a point $x \in X$ and choose $\lambda$ so that
    \[\Bun_N^{\lambda \cdot x} = \Bun_B \underset{\Bun_T}{\times} \{\omega^{\rhov}(\lambda \cdot x)\}\]
is a scheme. Write
    \[\WW_\lambda = \Sat_x(V_\lambda) \star \WW_\vac,\]
viewed a sheaf on $\Bun_N^{\lambda \cdot x}$. Since $\id \to V_\lambda^\vee \otimes V_\lambda$ is a direct summand, it suffices by $\Sph_{G,x}$-linearity to show that the cofiber of
    \[\Sat_x(V_\lambda) \star \eta_G \simeq \eta_{G,\Ran}(\WW_\lambda)\]
is cuspidal. In other words, we need to show for every proper parabolic $P$ that the morphism
\begin{equation}
\label{eqn:CT_eta}
    \CT_{P,!}^- \cdot \Poinc_{G,!} \cdot \WW_\lambda \to \CT_{P,!}^- \cdot \PsId_{G,!} \cdot \coeff_*^\vee \cdot \WW_{\lambda, \co}
\end{equation}
is an isomorphism.

\subsubsection{Third reduction} Let us rewrite the source of (\ref{eqn:CT_eta}) by use of the Jacquet functor
    \[
    \Jac_x : \Whit(G)_x \overset{\Jac_\Conf}{\longrightarrow} \Whit(M)_{\Conf_x} \overset{\sprd_*}{\longrightarrow} \Whit(M)_{\Ran_x}
    \]
and the coJacquet functor
    \[
    \co\Jac_x : \Whit(M)_{\Ran_x} \overset{\sprd^!}{\longrightarrow} \Whit(M)_{\Conf_x} \overset{\Jac_{\Conf,\co}}{\longrightarrow} \Whit(G)_x.
    \]
We obtain
\begin{equation}
\label{eqn:3_morphism_1}
\begin{split}
    & \CT_{P,!}^- \cdot \Poinc_{G,x,!} \cdot \WW_\lambda  \\
    & \qquad\simto \Poinc_{M,\Ran_x,!} \cdot \Jac_x \cdot \WW_\lambda \\
    & \qquad\qquad \simto \PsId_{M,!} \cdot \coeff_{M,\Ran_x,*}^\vee \cdot \left(\Jac_x \cdot \WW_\lambda \right)_\co \, [2\delta_M]\\
    & \qquad\qquad\qquad \simto \PsId_{M,!} \cdot \coeff_{M,\Ran_x,*}^\vee \cdot \co\Jac_x^\vee \cdot \WW_{\lambda, \co} \, [2\delta_G]
\end{split}
\end{equation}
The second arrow is induced the map $\eta_{M,\Ran_x}$, which is an isomorphism by induction. On the other hand, by rewriting the target of \ref{eqn:CT_eta}, we obtain
\begin{equation}
\label{eqn:3_morphism_2}
\begin{split}
    & \CT_{P,!}^- \cdot \Poinc_{G,x,!} \cdot \WW_\lambda \\
    & \qquad \to \CT_{P,!}^- \cdot \PsId_{G,!} \cdot \coeff_{G,x,*}^\vee \cdot \WW_{\lambda, \co} \, [2\delta_G] \\
    & \qquad\qquad \simto \PsId_{M,!} \cdot (\Eis_{P,!}^-)^\vee \cdot \coeff_{G,x,*}^\vee \cdot \WW_{\lambda, \co} \, [2\delta_G]\\
    & \qquad\qquad\qquad \simto \PsId_{M,!} \cdot \coeff_{M,\Ran_x, *}^\vee \cdot \Jac_{x,\co}^\vee \cdot \WW_{\lambda, \co} \, [2\delta_G].
\end{split}
\end{equation}
Therefore, it remains to perform the following:

\begin{construction}
\label{construction:3_morphism}
There is an identification of the morphism $(\ref{eqn:3_morphism_1})$ with $(\ref{eqn:3_morphism_2})$.
\end{construction}

\subsection{Construction of the identification}

\subsubsection{} We will compare (\ref{eqn:3_morphism_1}) with (\ref{eqn:3_morphism_2}) by describing both in explicit geometric terms. Consider the open substack
    \[
    \left( \Bun_N^{\lambda \cdot x} \underset{\Bun_G}{\times} \Bunt_{P^-} \underset{\Bun_M}{\times} \Bunc_M \right)^\trans
    \mon \Bun_N^{\lambda \cdot x} \underset{\Bun_G}{\times} \Bunt_{P^-} \underset{\Bun_M}{\times} \Bunc_M
    \]
cut out by the conditions:
\begin{itemize}
    \item the $N$ and $P^-$ structures are generically transverse, and
    \item the induced generic $N(M)$ structure is transverse to the $\Bunc_M$-structure.
\end{itemize}
This stack is equipped with with a natural projection $\pi$ onto $\Bun_M$.

\subsubsection{Description of (\ref{eqn:3_morphism_1})} Using the fact that $\WW_\lambda$ is a clean extension along the inclusion $\Bun_N^{\lambda \cdot x} \mon \Bun_N^{\pol\,x}$, it follows from a routine diagram chase that (\ref{eqn:3_morphism_1}) identifies with
    \[
    \pi_! \left( \WW_\lambda \stoximes j_{P,!}^- \stoximes r_{M,!} \right)^\trans 
    \to \pi_* \left( \WW_\lambda \shoximes j_{P,!}^- \shoximes r_{M,!} \right)^\trans [2 \delta_G + 2 \delta_M].
    \]
We shall identify (\ref{eqn:3_morphism_2b}) with this morphism.

\subsubsection{Functors defined by kernels} To handle (\ref{eqn:3_morphism_2}), we re-interpret the morphism
    \[
    \eta_G : \Poinc_{G,!} \cdot \WW_\lambda 
    \to \PsId_{G,!} \cdot \coeff_{x, *}^\vee \cdot \WW_{\lambda, \co} \, [2 \delta_G]
    \]
in terms of functors defined by kernels. Namely, we view
    \[
    \coeff_{\lambda, *} = \Gamma\left( \Bun_N^{\lambda \cdot x}, \WW_\lambda^{-1} \shimes \pi^!(-) \right)
    \qquad
    \coeff_{\lambda, !} = \Gamma_c\left(\Bun_N^{\lambda \cdot x}, \WW_\lambda^{-1} \stimes \pi^*(-) \right)
    \]
as functors defined by kernels
    \[
    \Sh(\Bun_G) \to \Sh(\pt).
    \]
By the assumption that $\Bun_N^{\lambda \cdot x}$ is a scheme, there is a canonical morphism
    \[\coeff_{\lambda, !} \to \coeff_{\lambda, *} \, [2 \delta_G].\]
Then $\eta_G$ identifies with the morphism
    \[
    (\id \boxtimes \coeff_{\lambda,!}) \cdot \Delta_{G,!} 
    \to (\id \boxtimes \coeff_{\lambda, *}) \cdot \Delta_{G,!} \, [2 \delta_G].
    \]

\subsubsection{Description of (\ref{eqn:3_morphism_2})} Now, factor (\ref{eqn:3_morphism_2}) as the composite of
\begin{equation}
\label{eqn:3_morphism_2a}
\begin{split}
    & \CT_{P,!}^- \cdot ( \Poinc_{G,!} \cdot \WW_\lambda ) \\
    & \qquad \to \CT_{P,!}^- \cdot ( \PsId_{G,!} \cdot \coeff_{x,*}^\vee \cdot \WW_{\lambda, \co} ) \, [2\delta_G] \\
    & \qquad\qquad \simto (\CT_{P,!}^- \cdot \PsId_{G,!}) \cdot \coeff_{x,*}^\vee \cdot \WW_{\lambda, \co} \, [2\delta_G] \\
    & \qquad\qquad\qquad \simto \PsId_{M,!} \cdot (\Eis_{P,!}^-)^\vee \cdot \coeff_{x,*}^\vee \cdot \WW_{\lambda, \co} \, [2\delta_G]
\end{split}
\end{equation}
and
\begin{equation}
\label{eqn:3_morphism_2b}
\begin{split}
    & \PsId_{M,!} \cdot (\Eis_{P,!}^-)^\vee \cdot \coeff_{x,*}^\vee \cdot \WW_{\lambda, \co} \, [2\delta_G] \\
    & \qquad \simto \PsId_{M,!} \cdot \coeff_{\Ran_x, *}^\vee \cdot \co\Jac_x^\vee \cdot \WW_{\lambda, \co} \, [2\delta_G]
\end{split}
\end{equation}
Then (\ref{eqn:3_morphism_2a}) identifies with the upper circuit of the diagram
\[
\begin{tikzcd}[
  column sep={18em,between origins},
  row sep={4em,between origins}
]
    (\CT_{P,!} \boxtimes \id) \cdot (\id \boxtimes \coeff_{\lambda,!}) \cdot \Delta_{G,!} 
    \arrow[r] \arrow[d] & (\CT_{P,!} \boxtimes \id) \cdot (\id \boxtimes \coeff_{\lambda, *}) \cdot \Delta_{G,!} \, [2 \delta_G] \arrow[d] \\
    (\id \boxtimes \coeff_{\lambda,!}) \cdot (\CT_{P,!} \boxtimes \id) \cdot \Delta_{G,!} 
    \arrow[r] \arrow[d] & (\id \boxtimes \coeff_{\lambda, *}) \cdot (\CT_{P,!} \boxtimes \id) \cdot \Delta_{G,!} \, [2 \delta_G] \arrow[d] \\
    (\id \boxtimes \coeff_{\lambda,!}) \cdot (\id \boxtimes \Eis_{P,!}^-) \cdot \Delta_{M,!}
    \arrow[r] & (\id \boxtimes \coeff_{\lambda, *}) \cdot (\id \boxtimes \Eis_{P,!}^-) \cdot \Delta_{M,!} \, [2 \delta_G]
\end{tikzcd}
\]
The lower face of this diagram is tautologically commutative. The commutation of the upper face is formal but tedious. Now reading the lower circuit of this diagram identifies (\ref{eqn:3_morphism_2a}) with
    \[
    \pi_! \left\{ \WW_\lambda \stoximes \left( j_{P,!}^- \stoximes r_{M,!} \right) \right\}^\trans 
    \to \pi_* \left\{ \WW_\lambda \stoximes \left( j_{P,!}^- \shoximes r_{M,!} \right) \right\}^\trans [2 \delta_G].
    \]
A further diagram chase identifies (\ref{eqn:3_morphism_2b}) with
    \[
    \pi_! \left\{ \WW_\lambda \stoximes \left( j_{P,!}^- \shoximes r_{M,!} \right) \right\}^\trans [2 \delta_G]
    \to \pi_* \left\{ \WW_\lambda \shoximes \left( j_{P,!}^- \shoximes r_{M,!} \right) \right\}^\trans [2 \delta_G + 2 \delta_M].
    \]
Composing these identifications completes Construction \ref{construction:3_morphism} and hence the proof of the main theorem.